\documentclass[11pt]{article}
\usepackage{amssymb, amsmath}
\usepackage{amscd}
\usepackage{supertabular}
\usepackage{setspace}

\usepackage{longtable}
\begin{document}
\newtheorem{Mthm}{Main Theorem.}
\newtheorem{Thm}{Theorem}[section]
\newtheorem{Prop}[Thm]{Proposition}
\newtheorem{Lem}[Thm]{Lemma}
\newtheorem{Cor}[Thm]{Corollary}
\newtheorem{Def}[Thm]{Definition}
\newtheorem{Guess}[Thm]{Conjecture}
\newtheorem{example}[Thm]{Example}
\newtheorem{Rmk}{Remark.}
\newtheorem{Not}{Notation.}
\newtheorem{thmA}{Theorem A.}
\newtheorem{thmB}{Theorem B.}
\newtheorem{thmC}{Theorem C.}
\newtheorem{DefProp}{Proposition-Definition}
\renewcommand{\theThm} {\thesection.\arabic{Thm}}
\renewcommand{\theProp}{\thesection.\arabic{Prop}}
\renewcommand{\theLem}{\thesection.\arabic{Lem}}
\renewcommand{\theCor}{\thesection.\arabic{Cor}}
\renewcommand{\theDef}{\thesection.\arabic{Def}}
\renewcommand{\theGuess}{\thesection.\arabic{Guess}}
\renewcommand{\theDefProp}{\thesection.\arabic{DefProp}}
\renewcommand{\theRmk}{}\renewcommand{\theMthm}{}
\renewcommand{\theNot}{}\renewcommand{\thethmA}{}
\renewcommand{\thethmB}{}\renewcommand{\thethmC}{}
\renewcommand{\thefootnote}{\fnsymbol{footnote}}
\newcommand{\spec}{\textnormal{Spec}\hspace{1mm}}
\newcommand{\proj}{\textnormal{Proj}\hspace{1mm}}
\newcommand{\Diff}{\textnormal{Diff}\hspace{1mm}}
\newcommand{\pf}{{\bfseries\itshape Proof. }}
\newcommand{\mult}{\textnormal{mult}}
\newcommand{\com}{\hspace{-2mm}\textnormal{\textbf{.}}\hspace{2mm}}
\newcommand{\bir}{-\hspace{-1mm}\rightarrow}
\newcommand{\qed}{\hfill $\Box$\newline}
\newcommand{\lct}{\operatorname{lct}}
\newcommand{\tlct}{\operatorname{tlct}}

\Large
\begin{center}
\textbf{Log canonical thresholds on Gorenstein canonical del Pezzo surfaces}
\end{center}
\vspace{4mm}

\normalsize
\begin{center}
\textbf{Jihun Park and Joonyeong Won
}\\
\vspace{5mm} \small{\itshape Department of Mathematics, POSTECH,
Pohang, Kyungbuk, $790$-$784$, Republic of Korea.}
\end{center}
\noindent
\hrulefill \vspace{2mm}

\begin{center}
\begin{minipage}[m]{.93\linewidth}
 \small {\emph{Abstract.}\ \ We classify all the effective anticanonical divisors on weak del Pezzo surfaces. Through this classification we obtain the smallest number among the log canonical thresholds of effective anticanonical divisors on a given Gorenstein canonical del Pezzo surface.
 \vspace{2mm}

  }

  \vspace{2mm}

 \emph{$2000$ Mathematics subject classification}: Primary 14J26; 14J17
\end{minipage}
\end{center}

\noindent
\hrulefill \normalsize

 \footnotetext[1]{Corresponding author : J. Park, Tel.
 +82-54-279-2059 Fax. +82-54-279-2799.}
 \footnotetext[2]{E-mail addresses :
\texttt{wlog@postech.ac.kr} (J. Park); \texttt{leonwon@postech.ac.kr} (J. Won)}
 \footnotetext[3]{This work has been supported by the Korea
Research Foundation Grant
(KRF-2007-412-J02302).}
\thispagestyle{empty}

\section{Introduction.}
Unless otherwise mentioned, all varieties are assumed to be projective, normal, and  defined
over~$\mathbb{C}$.

Let $X$ be a variety with at worst log canonical singularities    and $D$ be an
effective divisor on $X$.  The
log canonical threshold $c_x(D)$ of $D$ at a point $x$ in $X$ is defined by
$$c_x(D)=\mathrm{sup}\big\{c\ \big|\ \text{the log pair }(X,cD)\text{ is log canonical at the point $x$} \big\}.$$
The log canonical threshold $c(X,D)$ of the divisor $D$  is defined by
$$c(X,D)=\sup\big\{c\ \big| \ \text{the log pair }(X,cD)\text{ is log canonical} \big\}=\inf_{x\in X}\big\{ c_x(D) \big\}.$$
It is known that these numbers can be defined by some other equivalent ways. For instance, if the divisor $D$ is defined by a regular function $f$ near a smooth point $x$, then  the log canonical
threshold $c_x(D)$ of $D$ at the point $x$ is the number defined by
$$c_x(D)=\mathrm{sup}\big\{c\ \big|\ |f|^{-2c} \text{ is locally integrable near }x  \big\}.$$
It is also related to the Bernstein-Sato polynomial of the regular function $f$.
The log canonical threshold $c(X,D)$ of the divisor $D$  can be  defined by multiplier ideal sheaf as follows:
$$c(X,D)=\sup\big\{ c\ \big| \ \mathcal{J}(cD)=\mathcal{O}_X \big\},$$
where $\mathcal{J}(cD)$ is the multiplier ideal sheaf of $cD$.

Log canonical threshold, like multiplicity, measures how singular a divisor is.  However, log canonical threshold is a subtler
invariant than multiplicity. For instance, for the divisor~$D$ 
on $\mathbb{C}^2$  defined by $x^2=y^n$ around the origin $0$ that is of multiplicity $2$ at the origin, the log canonical threshold $c_0(D)$ of the divisor $D$ at the origin has different values as $n$
varies.
Also, log canonical threshold is  rather
difficult to calculate in a general case.  However, it has many
amazing properties and presents important applications to various areas such as
birational geometry and K\"{a}hler geometry.

The following theorem is one of the motivations of the present paper.

\begin{Thm}
\label{theorem:KE} Suppose that $X$ is a $n$-dimensional Fano orbifold. If
there is a positive real number $\epsilon$ such that for every effective $\mathbb{Q}$-divisor $D$ numerically equivalent to $-K_X$ the log pair~$(X, \frac{n+\epsilon}{n+1}D)$ is Kawamata log terminal, then
 $X$ has a K\"ahler--Einstein metric.
\end{Thm}
\pf See \cite{CSh08} and \cite{DeKo01}.\qed

This result  motivates the definition of the following numerical invariants.

\begin{Def}
Let $X$ be a Fano variety with at worst log terminal
singularities. The $m$-th global log canonical threshold of $X$ is defined by the number
$$
\mathrm{lct}_m\big(X\big)=\mathrm{sup}\left\{\lambda\in\mathbb{Q}\ \left|%
\aligned
&\text{the log pair}\ \left(X, \frac{\lambda}{m} D\right)\ \text{is log canonical }\\
&\text{for any effective divisor}\ D\in \Big|-mK_{X}\big|\\
\endaligned\right.\right\}.%
$$
The global log canonical threshold is defined by $\mathrm{lct}(X)=\inf\big\{ \ \mathrm{lct}_m(X)~ |~ m\in
\mathbb{N}\big\}$. Here, we do not define the $m$-th global log canonical threshold of $X$ if the linear system $|-mK_X|$ is empty.
\end{Def}
We can see that
$\mathrm{lct}(X)$ is the supremum of the  values $c$ such that the log
pair $(X,cD)$ is log canonical for every effective $\mathbb{Q}$-divisor
numerically equivalent to $-K_X$. Using the global log canonical threshold, Theorem~\ref{theorem:KE} can be reinterpreted as saying that
the Fano manifold $X$ admits a K\"ahler-Einstein metric if $$
\mathrm{lct}(X)>\frac{\mathrm{dim}(X)}{\mathrm{dim}(X)+1}.
$$

The paper \cite{Pu05} also shows that the global log canonical threshold plays important role in rationality problems.

The first global log canonical threshold may be a
cornerstone to get $\mathrm{lct}(X)$. It is natural that we  ask whether
there is an integer $m$ with $\mathrm{lct}_m(X)=\mathrm{lct}(X)$. We can find some evidence in simple cases.

\begin{Thm}\label{example:del Pezzo}
Let $X$ be a smooth del Pezzo surface. Then
\begin{equation*}
\mathrm{lct}\left(X\right)=\mathrm{lct}_{1}(X)=\left\{%
\aligned
&1/3\ \mathrm{ if}\ X\cong\mathbb{F}_{1}\ \mathrm{or}\ K_{X}^{2}\in\{7,9\},\\%
&1/2\ \mathrm{ if}\ X\cong\mathbb{P}^{1}\times\mathbb{P}^{1}\ \mathrm{or}\ K_{X}^{2}\in\{5,6\},\\%
&2/3\ \mathrm{if}\ K_{X}^{2}=4,\\%
&2/3\ \mathrm{if}\ X\ \mathrm{is\ a\ cubic\ in}\ \mathbb{P}^{3}\ \mathrm{with\ an\ Eckardt\ point},\\%
&3/4\ \mathrm{if}\ X\ \mathrm{is\ a\ cubic\ in}\ \mathbb{P}^{3}\ \mathrm{without\ Eckardt\ points},\\%
&3/4\ \mathrm{if}\ K_{X}^{2}=2\ \mathrm{and}\ |-K_{X}|\ \mathrm{has\ a\ tacnodal\ curve},\\%
&5/6\ \mathrm{if}\ K_{X}^{2}=2\ \mathrm{and}\ |-K_{X}|\ \mathrm{has\ no\ tacnodal\ curves},\\%
&5/6\ \mathrm{if}\ K_{X}^{2}=1\ \mathrm{and}\ |-K_{X}|\ \mathrm{has\ a\ cuspidal\ curve},\\%
&1\ \mathrm{if}\ K_{X}^{2}=1\ \mathrm{and}\ |-K_{X}|\ \mathrm{has\ no\ cuspidal\ curves}.\\%
\endaligned\right.%
\end{equation*}

\end{Thm}
\pf  See
\cite{Ch07a} and \cite{P99}. \qed

Throughout the present paper, an algebraic surface $S$ with ample
anticanonical divisor will be called a del Pezzo surface of
degree $d$ if it has at worst normal
Gorenstein canonical singularities and the self-intersection
number of the anticanonical divisor is $d$. Let $\Sigma$ be the set of singular points of $S$.
For singular del Pezzo
surfaces of degree 3, the paper \cite{Ch07b} shows the following:
\begin{Thm}\label{singular cubic}
 Suppose that $S$ is a cubic del Pezzo surface in
$\mathbb{P}^{3}$ and $\Sigma\ne\varnothing$. Then
$$\mathrm{lct}\big(S\big)=\mathrm{lct_1}\big(S\big)=\left\{%
\aligned & 1/6 \text{ if } \Sigma=\big\{\textnormal{E}_{6}\big\},\\%
& 1/4 \text{ if } \Sigma\supseteq\big\{\textnormal{A}_{5}\big\},\Sigma=\big\{\textnormal{D}_{5}\big\},\\%
& 1/3 \text{ if } \Sigma\supseteq\big\{\textnormal{A}_{4}\big\},\big\{2\textnormal{A}_{2}\big\},\Sigma=\big\{\textnormal{D}_{4}\big\},\\%
& 2/3 \text{ if }
\Sigma=\big\{\textnormal{A}_{1}\big\},\\
& 1/2 \text{ otherwise}.\\%
\endaligned\right.$$
\end{Thm}

For a del Pezzo surface $S$, the first global log canonical threshold $\mathrm{lct}_1(S)$ is meaningful by
itself. It has a nice application to birational maps between del Pezzo fibrations (see \cite{P99} or \cite{P01}).
The paper \cite{P01} has computed all the values of  $\mathrm{lct}_1(S)$ for  del Pezzo surfaces $S$ of degree 1.
\medskip

The aim of the present paper is to get all the values of the first global log canonical thresholds of Gorenstein canonical del Pezzo surfaces with singular points. This can be done by handling effective anticanonical divisors on the minimal resolutions of  del Pezzo surfaces.

To this end, we are first required to have information on
singularities of del Pezzo surfaces. This information can be
obtained from \cite{Dema80}, \cite{Furu86}, \cite{HiWa81},
\cite{MiZha88}, \cite{MiZha93}, \cite{U83}, and \cite{Zha88}.
Instead of studying singularities of del Pezzo surfaces, we are
able to understand them by studying configurations of $-2$-curves
on smooth surfaces with nef and big anticanonical divisor,
so-called weak del Pezzo surfaces. We can also obtain information
on effective anticanonical divisors on (weak) del Pezzo surfaces
in such a way.

We have the following geometric descriptions for del Pezzo surfaces which provide us with the relation to  weak del Pezzo surfaces
and a way to construct effective anticanonical divisors.

\begin{Thm}
Let $S$ be a del Pezzo surface of
degree $d$. Then
\begin{enumerate}
\item $1\leq d\leq 9$. \item If $d=9$, then
    $S\cong\mathbb{P}^2$.
\item If $d=8$, then either
    $S\cong\mathbb{P}^1\times\mathbb{P}^1$ or
    $S\cong\mathbb{F}_1$ or $S$ is a cone over
    a quadric in $\mathbb{P}^2$. \item If $1\leq d\leq 7$,
    then there exists a set of points in almost general
    position on $\mathbb{P}^2$ such that the blow up centered
    on the set is the minimal resolution of $S$.
\end{enumerate}
When we say that a finite set $\{p_1, p_2,\cdots,p_n\}$
 of points on the projective plane $\mathbb{P}^2$ (infinitely near points allowed)
is in almost general position, it means the following;
\begin{itemize}
\item No four of them are on a line.
\item No seven of them are on
a conic.
\item For all $j$ ($1\leq j\leq n-1$), the point
$p_{j+1}$ on the blow up $V_j$ of $\mathbb{P}^2$ centered at
$\{p_1,p_2,\cdots,p_j\}$ does not lie on any strict transform
$\hat{E}_i$ of $E_i$ ($1\leq i\leq j$) such that $\hat{E}_i^2=-2$,
where $E_i$ is an exceptional divisor on $V_i$.
\end{itemize}
\end{Thm}
\pf See \cite{Dema80} or \cite{HiWa81}. \qed

As we mentioned above, the singularities on a del Pezzo surface
can be described by the configuration of $-2$-curves on the
minimal resolution of the del Pezzo surface that is a  weak del
Pezzo surface. Meanwhile, the configurations of $-2$-curves can be
shown effectively  by their dual graphs (Dynkin diagrams).
\begin{Thm}
The singularities of a del Pezzo surface of degree $d$ are one of the following:
\begin{itemize}
\item[] $d=7$ ~~~~$\textnormal{A}_1$,
\item[] $d=6$ ~~~~any subgraph of
$\textnormal{A}_1+\textnormal{A}_2$,
\item[] $d=5$ ~~~~any proper subgraph of the extended Dynkin diagram
$\widetilde{\textnormal{A}}_4$,
\item[]  $d=4$ ~~~~any proper subgraph of the extended Dynkin diagram
$\widetilde{\textnormal{D}}_5$,
\item[] $d=3$ ~~~~any proper subgraph of the extended Dynkin diagram
$\widetilde{\textnormal{E}}_6$,
\item[] $d=2$ ~~~~$6\textnormal{A}_1$,
$\textnormal{D}_4+3\textnormal{A}_1$ or any proper subgraph of the
extended Dynkin diagram $\widetilde{\textnormal{E}}_7$.
\end{itemize}

\end{Thm}
\pf See \cite{HiWa81} and \cite{Pi77}.\qed

For the sake of the first global log canonical thresholds, we need
to distinguish some singularity types of del Pezzo surfaces of degree
$2$ with the same dual graphs.
To do so, we will distinguish $\textnormal{A}_5$ singularities into two types.
One has a
$-1$-curve intersecting the $-2$-curve corresponding to the vertex $v$ in the
dual graph of $\textnormal{A}_5$ such that
$\textnormal{A}_5-v=2\textnormal{A}_2$ on the minimal resolution of the del Pezzo surface. The other does not. In the former case the
type of singularities will be denoted by $\textnormal{A}_5'$
and in the latter
case by $\textnormal{A}_5''$. For singularity types $\textnormal{A}_5$  and
$\textnormal{A}_5+\textnormal{A}_1$  on del Pezzo surfaces of degree 2, there are two types for each
(see \cite{U83}). One is for $\textnormal{A}_5'$ and the other is for $\textnormal{A}_5''$. For singularity type
$\textnormal{A}_5+\textnormal{A}_2$ on del Pezzo surfaces of degree 2, there are only one type
(see \cite{U83}). The singularity $\textnormal{A}_5$ in this type is  $\textnormal{A}_5'$.

Also, there are two
types of singularities on del Pezzo surfaces of degree 2 with the
dual graph $3\textnormal{A}_1$ (resp. $4\textnormal{A}_1$) (see
\cite{U83}). One has a $-1$-curve on the del Pezzo surface which passes
through three $\textnormal{A}_1$ singular points (denoted by
$(3\textnormal{A}_{1})'$ (resp. $(4\textnormal{A}_{1})'$). The
other does not (denoted by $(3\textnormal{A}_{1})''$ (resp.
$(4\textnormal{A}_{1})''$). For singularity type
$\textnormal{A}_2+3\textnormal{A}_1$ on del Pezzo surfaces of degree~2, there are only one type
(see \cite{U83}). The singularities $3\textnormal{A}_1$ in this type are  $(3\textnormal{A}_1)'$.

We are now at the stage where we can state the main theorem of the present paper.
\begin{Thm}
\label{theorem:main} Let $S_d$ be a  del Pezzo surface of
degree $d$ and $\Sigma_d$ be the set of singular points in~$S_d$.
Suppose $\Sigma_d\neq \emptyset$. Then
\begin{itemize}
\item[]
$\mathrm{lct_1}\big(S_2\big)=\left\{%
\aligned & 1/6\text{ if }
\Sigma_2=\big\{\textnormal{E}_{7}\big\},\\
& 1/4 \text{ if } \Sigma_2=\big\{\textnormal{E}_{6}\big\},\Sigma_2\supseteq\big\{\textnormal{D}_{6}\big\},\\%
& 1/3 \text{ if }
\Sigma_2\supseteq\big\{\textnormal{D}_{5}\big\},\big\{\textnormal{A}_{5}'\big\},\\
& 1/2 \text{ if } \Sigma_2\supseteq\big\{(3\textnormal{A}_{1})'\big\},\big\{(4\textnormal{A}_{1})'\big\}
,\big\{5\textnormal{A}_{1}\big\}, \big\{\textnormal{A}_{3}\big\},\big\{\textnormal{A}_{4}\big\}, \big\{\textnormal{A}_{5}''\big\},
\big\{\textnormal{A}_{6}\big\}, \big\{\textnormal{A}_{7}\big\},\big\{\textnormal{D}_{4}\big\},
\\
& 2/3 \text{ otherwise. }\\%
\endaligned\right.$

\item[]
$\mathrm{lct_1}\big(S_3\big)=\left\{%
\aligned & 1/6 \text{ if } \Sigma_3=\big\{\textnormal{E}_{6}\big\},\\%
& 1/4 \text{ if } \Sigma_3\supseteq\big\{\textnormal{A}_{5}\big\},\Sigma_3=\big\{\textnormal{D}_{5}\big\},\\%
& 1/3 \text{ if } \Sigma_3\supseteq\big\{\textnormal{A}_{4}\big\},\big\{2\textnormal{A}_{2}\big\},\Sigma_3=\big\{\textnormal{D}_{4}\big\},\\%
& 2/3 \text{ if }
\Sigma_3=\big\{\textnormal{A}_{1}\big\},\\
& 1/2 \text{ otherwise}.\\%
\endaligned\right.$

\item[]
$\mathrm{lct_1}\big(S_4\big)=\left\{%
\aligned & 1/6 \text{ if } \Sigma_4=\big\{\textnormal{D}_{5}\big\},\\%
& 1/4 \text{ if } \Sigma_4\supseteq\big\{\textnormal{A}_{1}+\textnormal{A}_{3}\big\}, \Sigma_4=\big\{\textnormal{A}_{4}\big\},\Sigma_4=\big\{\textnormal{D}_{4}\big\},\\%
& 1/3 \text{ if } \Sigma_4=\big\{\textnormal{A}_{3}\big\},\Sigma_4\supseteq\big\{\textnormal{A}_{1}+\textnormal{A}_{2}\big\},\\%
& 1/2 \text{ otherwise.}\\%
\endaligned\right.$

\item[]
$\mathrm{lct_1}\big(S_5\big)=\left\{%
\aligned & 1/6 \text{ if } \Sigma_5=\big\{\textnormal{A}_{4}\big\},\\%
& 1/4 \text{ if }
\Sigma_5=\big\{\textnormal{A}_{3}\big\},\Sigma_5=\big\{\textnormal{A}_{1}+\textnormal{A}_{2}\big\},\\
& 1/3 \text{ if } \Sigma_5=\big\{\textnormal{A}_{2}\big\},\big\{2\textnormal{A}_{1}\big\},\phantom{\big\{\textnormal{A}_{1},\textnormal{A}_{1}\big\},}\\%
& 1/2 \text{ if } \Sigma_5=\big\{\textnormal{A}_{1}\big\}.\\%
\endaligned\right.$

\item[]
$\mathrm{lct_1}\big(S_6\big)=\left\{%
\aligned & 1/6 \text{ if } \Sigma_6=\big\{\textnormal{A}_{1}+\textnormal{A}_{2}\big\}.\\%
& 1/4 \text{ if }
\Sigma_6=\big\{\textnormal{A}_{2}\big\},\Sigma_6=\big\{2\textnormal{A}_{1}\big\},\\%
& 1/3 \text{ if }
\Sigma_6=\big\{\textnormal{A}_{1}\big\}.\\
\endaligned\right.$

\item[]
$\mathrm{lct_1}\big(S_7\big)=~~~1/4 \text{ if }
\Sigma_7=\big\{\textnormal{A}_{1}\big\}.$

\end{itemize}
\end{Thm}

The first log canonical thresholds of del Pezzo surfaces of degree
$1$ have been dealt with in~\cite{P01}.

Let $\pi :\tilde{S}\rightarrow S$ be the minimal resolution of
$S$.  Since we assume that the del Pezzo surface admits only Gorenstein canonical singularities, the resolution
 $\pi$ is crepant, i.e., $K_{\widetilde{S}}=\pi^*(K_S).$ Therefore the pull-back
 $\tilde{D}:=\pi^*(D)$ of an effective anticanonical divisor $D$ on $S$
 is an effective anticanonical divisor on $\tilde{S}$.
 We can write $\pi^*(D)=\bar{D}+E$,
where $\bar{D}$ is the strict transform of $D$ and
$\mathrm{Supp}(E)$ consists of $-2$-curves.

\begin{Lem}\label{counting}
If $D=\sum a_iD_i$ is an effective anticanonical divisor on a
del Pezzo surface of degree~$d$, then $\sum a_i\leq d$.
\end{Lem}
\pf It is easy to check. \qed

For a constant $c$,
$$\pi^*(K_S+cD)=K_{\tilde{S}}+c\tilde{D}.$$ Thus it is sufficient
to consider $\tilde{D}$ on $\tilde{S}$ to compute
$\mathrm{lct}_1(S)$. An effective anticanonical divisor $D$ on $S$
which  does not pass through any singular point of $S$ is not
different from the pull-back of~$D$ via $\pi$.  When we consider
effective anticanonical divisors which pass through singular
points of $S$, it suffices to investigate effective anticanonical
divisors  on the
 weak del Pezzo surface $\tilde{S}$ that contains a $-2$-curve.

Since  $\mathrm{lct}_1(S)$ is always at most $1$,  classifying all effective anticanonical divisors in $\tilde{S}$
that has either at least one component with multiplicity $\geq 2$ or components that are not normal crossing,
we can prove the main theorem. The classification will
be presented in the following section.

\section{The configuration of the anticanonical divisors.}\label{configuration}

Let us  summarize  the results of \cite{P01}. Let $S$ be a
 del Pezzo surface of degree $1$ and $\pi:\tilde{S}\to
S$ be the minimal resolution of $S$. Then for an effective
anticanonical divisor $D$ on $S$, the configuration of its
pull-back divisor $\tilde{D}$ by the morphism $\pi$ coincides with
one of Kodaira's elliptic fibers. Lemma~\ref{counting}
implies that every effective anticanonical divisor on $S$ is
irreducible and reduced since the degree of $S$ is $1$.  If the
divisor $D$ passes through a singular point of $S$ and the divisor~$\tilde{D}$ has a multiple component,
 then we can see the dual graph of the divisor $\tilde{D}$  must be
 one of those that appear in the tables  in
 Propositions~\ref{proposition:wt6},~\ref{proposition:wt4},~\ref{proposition:wt3},
 and~\ref{proposition:wt2}.
 For del Pezzo surfaces of degree $1$ we can
obtain all the configurations of $\tilde{D}$ from \cite{P01}.

Every effective anticanonical divisor on a weak del Pezzo surface of degree $d$ can be obtained
from an effective anticanonical divisor on a weak del Pezzo surface of degree $d+1$ via a suitable blow-up.

Suppose that we have obtained the list of all the configurations of effective anticanonical divisors on weak del Pezzo surfaces of degree $d$.
 Let $D_{d+1}= \sum a_iE_i+\sum b_jD_j\sim -K_{d+1}$ be an effective anticanonical divisor of a weak del Pezzo surface
 of degree $d+1$, where each $E_i$ is a $-2$-curve and $D_j$ is a prime divisor that is not a $-2$-curve.  In addition, let $\hat{D}=\sum D_j$. Then we consider the blow up
 $\psi$ of the weak del Pezzo surface of degree $d+1$ at a point $p$ in $D_k\setminus\cup E_i$  with $\mult_p(\hat{D})=1$. Then it produces
a new effective anticanonical divisor on the weak del Pezzo surface of degree $d$ that is the blow up by $\psi$. To be precise, we obtain
$$
\psi^*(-K_{d+1})-F=\sum a_i\bar{E}_i+\sum
b_j\bar{D}_j+(b_k-1)F\sim -K_{d},$$ where $\bar{E}_i$, $\bar{D}_j$
are the strict transforms of $E_i$, $D_j$, respectively, and  $F$
is the exceptional divisor of $\psi$. This effective anticanonical
divisor should appear in the list of all the configurations of
effective anticanonical divisors on weak del Pezzo surfaces of
degree $d$ that we have already obtained. Therefore, in order to
obtain the list of all the configurations of effective
anticanonical divisors on weak del Pezzo surfaces of degree $d+1$
from the list for degree $d$, we have only to consider two cases,
depending on $b_k=1$ or $b_k>1$, as follows:
\begin{itemize}
\item[\textbf{P-1.}] Suppose that we have an effective
anticanonical divisor on a weak del Pezzo surface of degree $d$
whose dual graph $\Gamma_d$ has a vertex  $v$  with weight 1 and
self-intersection number $j\geq -2$. Then there is the possibility
that  we have an effective anticanonical divisor on a weak del
Pezzo surface of degree $d+1$ whose dual graph  is the same as
$\Gamma_{d}$ except that the vertex $v$ has  self-intersection
number $j+1$. This is the case where $b_k=1$ in the description above.
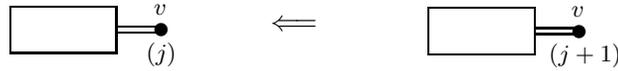
\begin{center}
\begin{minipage}[m]{.30\linewidth}\setlength{\unitlength}{.20mm}
\begin{center}
\begin{picture}(50,75)(120,12)
\put(130,30){\framebox(70,30)}\put(200,47){\line(1,0){30}}\put(200,43){\line(1,0){30}}\put(230,45){\circle*{8}}
 \put(220,25){\mbox{\footnotesize $(j)$}}
\put(225,55){\mbox{\footnotesize $v$}}
\end{picture}
\end{center}
\end{minipage}
\ \ \ \ \ \  $\Longleftarrow$ \ \ \ \ \ \ \ \ \
\begin{minipage}[m]{.30\linewidth}\setlength{\unitlength}{.20mm}
\begin{center}
\begin{picture}(500,75)(120,12)
\put(130,30){\framebox(70,30)}\put(200,47){\line(1,0){30}}\put(200,43){\line(1,0){30}}\put(230,45){\circle*{8}}
 \put(210,25){\mbox{\footnotesize $(j+1)$}}\put(225,55){\mbox{\footnotesize $v$}}
\end{picture}
\end{center}

\end{minipage}
\end{center}

\item[\textbf{P-2.}] Suppose that we have an effective
anticanonical divisor on a weak del Pezzo surface of degree $d$
whose dual graph $\Gamma_d$ has a vertex  $w$  with weight $b\geq
1$ and self-intersection number $-1$ that has only one adjacent
vertex $v$. We also suppose that the vertex $v$ has weight
$b+1$ and self-intersection number $j\geq -2$. Then there is the
possibility that  we have an effective anticanonical divisor on a
weak del Pezzo surface of degree $d+1$ whose dual graph  is the
same as $\Gamma_{d}$ except that the vertex  $v$ has
self-intersection number $j+1$ and it has no $w$. This is the case where $b_k>1$ in the description above.

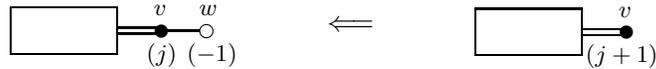
\begin{center}
\begin{minipage}[m]{.40\linewidth}\setlength{\unitlength}{.20mm}
\begin{center}
\begin{picture}(50,75)(120,12)
\put(130,30){\framebox(70,30)}\put(200,47){\line(1,0){30}}\put(200,43){\line(1,0){30}}\put(230,45){\circle*{8}}
 \put(220,25){\mbox{\footnotesize $(j)$}}\put(247,25){\mbox{\footnotesize $(-1)$}}
\put(225,55){\mbox{\footnotesize $v$}}\put(255,55){\mbox{\footnotesize $w$}}
\put(260,45){\circle{8}}\put(230,45){\line(1,0){25}}
\end{picture}
\end{center}

\end{minipage}
\ \ \ \ \ \ $\Longleftarrow$ \ \ \ \ \ \ \ \
\begin{minipage}[m]{.40\linewidth}\setlength{\unitlength}{.20mm}
\begin{center}
\begin{picture}(500,75)(120,12)
\put(130,30){\framebox(70,30)}\put(200,47){\line(1,0){30}}\put(200,43){\line(1,0){30}}\put(230,45){\circle*{8}}
 \put(203,25){\mbox{\footnotesize $(j+1)$}}
\put(225,55){\mbox{\footnotesize $v$}}
\end{picture}
\end{center}

\end{minipage}

\end{center}
\end{itemize}
Here the double lines mean that the vertex $v$ can be connected to either one vertex or more than one vertex.

Let $\tilde{S}$ be a weak del Pezzo surface of degree $d\leq 7$. Then there are $9-d$ points $p_{i,j}$, $i,j\geq 1$, on the projective plane $\mathbb{P}^2$ (infinitely near points allowed) in almost general position such that the blow up centered at these points is the surface $\tilde{S}$. We have a birational morphism $\pi:\tilde{S}\to\mathbb{P}^2$ that is a composition of a sequence of blow ups. Here each $p_{i,1}$ is a point on $\mathbb{P}^2$ and the point~$p_{i,j+1}$ is a point on the exceptional divisor of the blow up at the point $p_{i, j}$. The exceptional divisor of the blow up at the point $p_{i,j}$ is denoted by $E_{i,j}$. Then, we see
\[K_{\tilde{S}}=\pi(K_{\mathbb{P}^2})+\sum_{ij}jE_{i,j}.\]
For a divisor $D$ on $\mathbb{P}^2$, we have
\[\pi^*(D)=\bar{D}+\sum_{ij}(\mult_{p_{i,j}}(D))E_{i,j},\]
where $\bar{D}$ is the strict transform of $D$ by $\pi$. In
particular, a cubic curve $C=\sum a_hC_h$ (not necessarily
irreducible nor reduced) on $\mathbb{P}^2$ defines an
anticanonical divisor
\begin{equation}\label{adjuction}
\sum a_h\bar{C}_h+\sum_{i, j}\left(\sum_{k=1}^j\left(\mult_{p_{i,k}}(C)-1\right)\right)
E_{i,j}
\end{equation}
on the surface $\tilde{S}$.

For all the dual graphs in what follows, we will use the following
notation to distinguish smooth rational curves with various
self-intersection numbers;
\begin{center}
\begin{longtable}{|c|c|}
\hline
\begin{minipage}[m]{.25\linewidth}\setlength{\unitlength}{.20mm}
\begin{center}
\begin{picture}(10,10)(50,10)
\put(50,15){\circle*{5}}
\end{picture}
\end{center}
\end{minipage}&
\begin{minipage}[m]{.4\linewidth}
\begin{center}
$-2$-curve
\end{center}
\end{minipage}\\
\hline
\begin{minipage}[m]{.25\linewidth}\setlength{\unitlength}{.20mm}
\begin{center}
\begin{picture}(10,10)(50,10)
\put(50,15){\circle{5}}
\end{picture}
\end{center}

\end{minipage}&
\begin{minipage}[m]{.4\linewidth}
\begin{center}
$-1$-curve

\end{center}
\end{minipage}\\
\hline
\begin{minipage}[m]{.25\linewidth}\setlength{\unitlength}{.20mm}
\begin{center}
\begin{picture}(10,10)(50,10)
\put(48,13){\framebox(4,4)}
\end{picture}
\end{center}

\end{minipage}&
\begin{minipage}[m]{.4\linewidth}
\begin{center}
$0$-curve
\end{center}
\end{minipage}\\
\hline
\begin{minipage}[m]{.25\linewidth}\setlength{\unitlength}{.20mm}
\begin{center}
\begin{picture}(10,10)(50,10)
\put(50,15){\circle*{4}} \put(48,13){\framebox(4,4)}
\end{picture}
\end{center}

\end{minipage}&
\begin{minipage}[m]{.4\linewidth}
\begin{center}
$1$-curve
\end{center}
\end{minipage}\\
\hline
\begin{minipage}[m]{.25\linewidth}\setlength{\unitlength}{.20mm}
\begin{center}
\begin{picture}(10,10)(50,10)
\put(50,15){\circle{4}} \put(48,13){\framebox(4,4)}
\end{picture}
\end{center}

\end{minipage}&
\begin{minipage}[m]{.4\linewidth}
\begin{center}
$2$-curve
\end{center}
\end{minipage}\\

\hline
\end{longtable}
\end{center}
The number near each vertex is the multiplicity of the curve
corresponding to the vertex. The number $1$ for multiplicity $1$
will be always omitted.

In the tables of
Propositions~\ref{proposition:wt6},~\ref{proposition:wt4},~\ref{proposition:wt3},
 and~\ref{proposition:wt2}, some rows are marked with~$\surd$. We shall see in Section~3 that only those entries marked with ~$\surd$ need be considered in order to prove Theorem~\ref{theorem:main} (see Propositions~\ref{change}, \ref{negative-curves-only} and~\ref{no-1}).

The columns labeled as Example carry configurations of
divisors on certain blow ups of $\mathbb{P}^2$ in order to show
existence of anticanonical divisors of given types on  weak del
Pezzo surfaces. In each configuration, solid lines, which denote exceptional curves of blow-ups of $\mathbb{P}^2$, will show
the way of blow-ups of $\mathbb{P}^2$. Among the solid lines, thin
lines (always drawn horizontally) denote $-1$-curves and thick lines (always drawn slantingly) denote $-2$-curves. The
dotted curves in each configuration are the strict transformation
of  a cubic curve (not necessarily reduced nor irreducible) via
the blow-ups. The letters $L$ and $Q$ right beside the dotted
curves mean that each corresponding dotted curve is  the strict transformation of a line and an irreducible conic,
respectively. In addition, $2L$ and $3L$ mean the strict transformation of a double line and a
triple line, respectively.

Since each curve in an effective anticanonical divisor
on a weak del Pezzo surface of degree $1$ has weight at most $6$,
so does a curve in an effective  anticanonical divisor on a weak del Pezzo surface.

\begin{Prop}\label{proposition:wt6}
If the dual graph for an effective anticanonical divisor  on a weak del Pezzo surface  has a vertex of weight $6$,
then it is exactly one of the following:
\footnotesize

\begin{center}
\begin{longtable}{|l|c|c|}
\hline \multicolumn{3}{|c|}
{Degree $1$} \\
\hline

\phantom{$\surd$\ $2\textnormal{A}_1+\textnormal{A}_2$} &Configuration& Example \\
\hline

$\surd$ \ $\textnormal{E}_8$ &
\begin{minipage}[m]{.28\linewidth}\setlength{\unitlength}{.20mm}
\begin{center}
\begin{picture}(500,60)(65,10)
\multiput(70,40)(30,0){7}{\circle*{5}} \put(130,20){\circle*{5}}
\put(100,40){\line(1,0){177.5}} \put(130,40){\line(0,-1){20}}
\put(70,40){\line(1,0){30}}
\put(280,40){\circle{5}} \put(70,45){$2$} \put(100,45){$4$}
\put(130,45){$6$} \put(160,45){$5$} \put(190,45){$4$}
\put(220,45){$3$} \put(120,20){$3$} \put(250,45){$2$}
\end{picture}
\end{center}
\end{minipage}
&
\begin{minipage}[m]{.25\linewidth}\setlength{\unitlength}{.20mm}
\begin{center}
\begin{picture}(400,50)(-10,15)

\qbezier[15](45,42)(58,36)(71,30)
\put(41,45){\mbox{\footnotesize $3L$}}

\put(133,42){\line(1,0){30}}

\thicklines
\put(0,30){\line(2,1){30}}
\put(20,45){\line(2,-1){30}}
\put(40,30){\line(2,1){30}}
\put(60,45){\line(2,-1){30}}
\put(80,30){\line(2,1){30}}
\put(100,45){\line(2,-1){30}}
\put(120,30){\line(2,1){30}}

\end{picture}
\end{center}
\end{minipage}
\\
\hline
\end{longtable}
\end{center}

\begin{center}
\begin{longtable}{|l|c|c|}
\hline \multicolumn{3}{|c|}
{Degree $2$} \\
\hline

\phantom{$\surd$\ $2\textnormal{A}_1+\textnormal{A}_2$} &Configuration& Example \\
\hline

$\surd$ \ $\textnormal{E}_7$ &
\begin{minipage}[m]{.28\linewidth}\setlength{\unitlength}{.20mm}

\begin{center}
\begin{picture}(500,60)(90,10)
\multiput(100,40)(30,0){6}{\circle*{5}} \put(280,40){\circle{5}}
\put(160,20){\circle*{5}} \put(100,40){\line(1,0){177.5}}
\put(160,40){\line(0,-1){20}} \put(100,45){$2$} \put(130,45){$4$}
\put(160,45){$6$} \put(190,45){$5$} \put(220,45){$4$}
\put(250,45){$3$} \put(165,20){$3$} \put(280,45){$2$}
\end{picture}
\end{center}
\end{minipage}
&
\begin{minipage}[m]{.25\linewidth}\setlength{\unitlength}{.20mm}
\begin{center}
\begin{picture}(400,50)(-10,15)

\qbezier[15](45,42)(58,36)(71,30)
\put(41,45){\mbox{\footnotesize $3L$}}

\put(117,33){\line(1,0){30}}

\thicklines
\put(0,30){\line(2,1){30}}
\put(20,45){\line(2,-1){30}}
\put(40,30){\line(2,1){30}}
\put(60,45){\line(2,-1){30}}
\put(80,30){\line(2,1){30}}
\put(100,45){\line(2,-1){30}}

\end{picture}

\end{center}

\end{minipage}
\\
\hline
\end{longtable}
\end{center}

\begin{center}
\begin{longtable}{|l|c|c|}
\hline \multicolumn{3}{|c|}
{Degree $3$} \\
\hline

\phantom{$\surd$\ $2\textnormal{A}_1+\textnormal{A}_2$} &Configuration& Example \\
\hline

$\surd$ \ $\textnormal{E}_6$ &
\begin{minipage}[m]{.28\linewidth}\setlength{\unitlength}{.20mm}

\begin{center}
\begin{picture}(500,50)(90,8)
\multiput(100,40)(30,0){5}{\circle*{5}} \put(250,40){\circle{5}}
\put(160,20){\circle*{5}} \put(100,40){\line(1,0){147.5}}
\put(160,40){\line(0,-1){20}} \put(100,45){$2$} \put(130,45){$4$}
\put(160,45){$6$} \put(190,45){$5$} \put(220,45){$4$}
\put(250,45){$3$} \put(165,20){$3$}
\end{picture}
\end{center}
\end{minipage}
&
\begin{minipage}[m]{.25\linewidth}\setlength{\unitlength}{.20mm}
\begin{center}
\begin{picture}(400,50)(-10,15)

\qbezier[15](45,42)(58,36)(71,30)
\put(41,45){\mbox{\footnotesize $3L$}}

\put(97,43){\line(1,0){30}}

\thicklines
\put(0,30){\line(2,1){30}}
\put(20,45){\line(2,-1){30}}
\put(40,30){\line(2,1){30}}
\put(60,45){\line(2,-1){30}}
\put(80,30){\line(2,1){30}}

\end{picture}

\end{center}

\end{minipage}
\\
\hline
\end{longtable}
\end{center}

\begin{center}
\begin{longtable}{|l|c|c|}
\hline \multicolumn{3}{|c|}
{Degree $4$} \\
\hline

\phantom{$\surd$\ $2\textnormal{A}_1+\textnormal{A}_2$} &Configuration& Example \\
\hline

$\surd$ \ $\textnormal{D}_5$&
\begin{minipage}[m]{.28\linewidth}\setlength{\unitlength}{.20mm}
\begin{center}
\begin{picture}(500,60)(90,5)
\multiput(100,40)(30,0){4}{\circle*{5}} \put(220,40){\circle{5}}
\put(160,17.5){\circle*{5}} \put(100,40){\line(1,0){117.5}}
\put(160,40){\line(0,-1){20}} \put(100,45){$2$} \put(130,45){$4$}
\put(160,45){$6$} \put(190,45){$5$} \put(220,45){$4$}
\put(165,20){$3$}
\end{picture}
\end{center}
\end{minipage}

&
\begin{minipage}[m]{.25\linewidth}\setlength{\unitlength}{.20mm}
\begin{center}
\begin{picture}(400,50)(-10,10)
\qbezier[15](40,48)(55,38)(70,28)
\put(43,48){\mbox{\footnotesize $3L$}}

\put(78,32){\line(1,0){30}}
\thicklines
\put(0,30){\line(2,1){30}}
\put(20,45){\line(2,-1){30}}
\put(40,30){\line(2,1){30}}
\put(60,45){\line(2,-1){30}}

\end{picture}
\end{center}

\end{minipage}
\\
\hline
\end{longtable}
\end{center}

\begin{center}
\begin{longtable}{|l|c|c|}
\hline \multicolumn{3}{|c|}
{Degree $5$} \\
\hline

\phantom{$\surd$ \ $2\textnormal{A}_1+\textnormal{A}_2$} &Configuration& Example \\
\hline

$\surd$ \ $\textnormal{A}_4$&
\begin{minipage}[m]{.28\linewidth}\setlength{\unitlength}{.20mm}
\begin{center}
\begin{picture}(500,60)(90,5)
\multiput(100,40)(30,0){3}{\circle*{5}} \put(190,40){\circle{5}}
\put(160,17.5){\circle*{5}} \put(100,40){\line(1,0){87.5}}
\put(160,40){\line(0,-1){20}} \put(100,45){$2$} \put(130,45){$4$}
\put(160,45){$6$} \put(190,45){$5$} \put(165,20){$3$}
\end{picture}
\end{center}
\end{minipage}

&
\begin{minipage}[m]{.25\linewidth}\setlength{\unitlength}{.20mm}
\begin{center}
\begin{picture}(400,50)(-10,10)

\qbezier[15](40,48)(55,38)(70,28)
\put(43,48){\mbox{\footnotesize $3L$}}

\put(55,43){\line(1,0){30}}
\thicklines
\put(0,30){\line(2,1){30}}
\put(20,45){\line(2,-1){30}}
\put(40,30){\line(2,1){30}}

\end{picture}
\end{center}

\end{minipage}

\\
\hline
\end{longtable}
\end{center}

\begin{center}
\begin{longtable}{|l|c|c|}
\hline \multicolumn{3}{|c|}
{Degree $6$} \\
\hline

\phantom{$\surd$ \ $2\textnormal{A}_1+\textnormal{A}_2$} &Configuration& Example \\
\hline

$\surd$ \ $\textnormal{A}_2+\textnormal{A}_1$&
\begin{minipage}[m]{.28\linewidth}\setlength{\unitlength}{.20mm}
\begin{center}
\begin{picture}(500,60)(120,5)
\multiput(190,40)(30,0){1}{\circle{5}} \put(220,40){\circle*{5}}
 \put(127.5,40){\circle*{5}}\put(160,40){\circle*{5}}
\put(130,40){\line(1,0){27.5}}\put(162.5,40){\line(1,0){25}}\put(192.5,40){\line(1,0){25}}
 \put(130,45){$2$} \put(160,45){$4$}
 \put(190,45){$6$}
\put(220,45){$3$}
\end{picture}
\end{center}
\end{minipage}

&
\begin{minipage}[m]{.25\linewidth}\setlength{\unitlength}{.20mm}
\begin{center}
\begin{picture}(400,50)(-10,10)

\qbezier[15](40,48)(55,38)(60,25)
\put(43,48){\mbox{\footnotesize $3L$}}

\put(38,33){\line(1,0){30}}
\thicklines
\put(0,30){\line(2,1){30}}
\put(20,45){\line(2,-1){30}}

\end{picture}
\end{center}

\end{minipage}
\\
\hline
\end{longtable}
\end{center}

\normalsize
\end{Prop}
\pf
Since there is no $\textnormal{E}_8$ on a del Pezzo surface of degree $2$,
we have only \textbf{P-2} possibilities.
Furthermore, such a dual graph on  a del Pezzo surface of degree $\geq 2$ has no vertex with weight~1. Therefore, we also have only \textbf{P-2} possibilities after degree 2.
\qed

The weight $5$ never appears as a maximum weight in an effective
anticanonical divisor on a weak del Pezzo surface of degree $1$
and the maximum weight is preserved under the changes by
\textbf{P-1} and \textbf{P-2}. Therefore, if an effective
anticanonical divisor on a weak del Pezzo surface has multiplicity
$5$ along a curve, then it must have another curve along which it
has multiplicity $6$.

\begin{Lem}\label{wt4}
 If the dual graph for an effective anticanonical divisor  on a weak del Pezzo surface  has a vertex of weight $4$ as a maximum,
then it  contains at most one vertex with nonnegative self-intersection number. In such a case, the vertex has self-intersection number $0$.
\end{Lem}
\pf Let $\tilde{S}$ be a weak del Pezzo surface of degree $d\leq 7$.
Then it is obtained by suitable blow ups $\pi:\tilde{S}\to \mathbb{P}^2$.
Every effective anticanonical divisor on $\tilde{S}$ can be
obtained from a cubic curve $C=\sum a_hC_h$
(not necessarily irreducible nor reduced) on $\mathbb{P}^2$ as formula~(\ref{adjuction}) shows.
Then for some $(i, j)$, we must have $
\sum_{k=1}^j\left(\mult_{p_{i,k}}(C)-1\right)=4 $. Since $\mult_{p_{i,k}}(C)$ is non-increasing
as $k$ grows and the curve $C$ is cubic,  the possible sequences for
$\{\mult_{p_{i, k}}(C)\}_{k=1}^{j}$ have the form as follows:

$$(2,2,2,2, *, *,\cdots), \ (3,2,2, *, *,\cdots),\ (3,3, *, *,\cdots).$$
Therefore, the curve $C$ consists of only lines.

Since at most three $p_{i,k}$ can be on a
line, the first form is impossible, the second must be~$(3,2,2)$ and the last must be $(3,3)$.

For the sequence $(3,3)$, the curve $C$ is a triple line. Exactly two points of $p_{i,k}$ are over $C$ and hence the effective anticanonical divisor on $\tilde{S}$ given by $C$ has no nonnegative self-intersection curve.

For the sequence $(3,2,2)$, the curve $C$ consists of one double line and one single line.
We should take one blow up at the intersection point of the double line and
the single line. Then two more blow ups at some points over the double line must follow.
 Therefore, the effective anticanonical divisor given by $C$ has only one  non-negative
  self-intersection curve and its self-intersection number is $0$.
\qed

\begin{Prop}\label{proposition:wt4}
If the dual graph for an effective anticanonical divisor  on a weak del Pezzo surface  has a vertex of weight $4$ as a maximum,
then it is exactly one of the following:
\footnotesize

\begin{center}
\begin{longtable}{|l|c|c|}
\hline \multicolumn{3}{|c|}
{Degree $1$} \\
\hline

\phantom{$\surd$\ $2\textnormal{A}_1+\textnormal{A}_2$} &Configuration& Example \\
\hline

$\surd$\ $\textnormal{E}_7$ &
\begin{minipage}[m]{.25\linewidth}\setlength{\unitlength}{.20mm}

\begin{center}
\begin{picture}(500,60)(65,10)
\multiput(100,40)(30,0){6}{\circle*{5}} \put(160,20){\circle*{5}}
\put(100,40){\line(1,0){150}}
 \put(160,40){\line(0,-1){20}}
\put(72.5,40){\line(1,0){27.5}}
\put(70,40){\circle{5}}

\put(100,45){$2$} \put(130,45){$3$} \put(160,45){$4$}
 \put(190,45){$3$}
\put(220,45){$2$}
 \put(150,20){$2$}
\end{picture}
\end{center}
\end{minipage}
&
\begin{minipage}[m]{.25\linewidth}\setlength{\unitlength}{.20mm}
\begin{center}
\begin{picture}(400,50)(-10,10)

\qbezier[15](45,42)(58,36)(71,30)
\put(41,45){\mbox{\footnotesize $2L$}}
\qbezier[15](10,42)(10,28)(10,15)
\put(2,45){\mbox{\footnotesize $L$}}

\put(117,33){\line(1,0){30}}
\put(0,20){\line(1,0){30}}

\thicklines
\put(0,30){\line(2,1){30}}
\put(20,45){\line(2,-1){30}}
\put(40,30){\line(2,1){30}}
\put(60,45){\line(2,-1){30}}
\put(80,30){\line(2,1){30}}
\put(100,45){\line(2,-1){30}}

\end{picture}
\end{center}
\end{minipage}
\\
\hline
\end{longtable}
\end{center}

\begin{center}
\begin{longtable}{|l|c|c|}
\hline \multicolumn{3}{|c|}
{Degree $2$} \\
\hline

\phantom{$\surd$\ $2\textnormal{A}_1+\textnormal{A}_2$}&Configuration& Example \\

\hline

$\surd$\ $\textnormal{E}_6$ &
\begin{minipage}[m]{.25\linewidth}\setlength{\unitlength}{.20mm}
\begin{center}
\begin{picture}(500,60)(90,0)
\multiput(130,40)(30,0){5}{\circle*{5}} \put(190,20){\circle*{5}}
\put(97.5,40){\circle{5}} \put(282.5,40){\circle{5}}
\put(100,40){\line(1,0){180}}
 \put(190,40){\line(0,-1){20}}
\put(130,45){$2$}

\put(160,45){$3$} \put(190,45){$4$}
 \put(220,45){$3$}
\put(250,45){$2$}
\put(185,5){$2$}
\end{picture}
\end{center}

\end{minipage}&
\begin{minipage}[m]{.25\linewidth}\setlength{\unitlength}{.20mm}
\begin{center}
\begin{picture}(400,40)(-10,12)

\qbezier[15](45,42)(55,36)(65,30)
\put(41,45){\mbox{\footnotesize $2L$}}
\qbezier[15](5,45)(5,30)(5,15)
\put(7,45){\mbox{\footnotesize $L$}}

\put(0,20){\line(1,0){30}}
\put(98,42){\line(1,0){30}}

\thicklines
\put(0,30){\line(2,1){30}}
\put(20,45){\line(2,-1){30}}
\put(40,30){\line(2,1){30}}
\put(60,45){\line(2,-1){30}}
\put(80,30){\line(2,1){30}}

\end{picture}
\end{center}

\end{minipage}\\

\hline

$\surd$\ $\textnormal{D}_6$ &
\begin{minipage}[m]{.25\linewidth}\setlength{\unitlength}{.20mm}

\begin{center}
\begin{picture}(500,60)(90,10)
\multiput(130,40)(30,0){5}{\circle*{5}} \put(160,20){\circle*{5}}
\put(97.5,40){\circle{5}} \put(100,40){\line(1,0){150}}
 \put(160,40){\line(0,-1){20}}
\put(100,45){$2$} \put(130,45){$3$} \put(160,45){$4$}
 \put(190,45){$3$}
\put(220,45){$2$}
\put(165,20){$2$}
\end{picture}
\end{center}
\end{minipage}&
\begin{minipage}[m]{.25\linewidth}\setlength{\unitlength}{.20mm}
\begin{center}
\begin{picture}(400,60)(-10,-3)
\qbezier[30](45,42)(25,24)(5,3)
\put(41,45){\mbox{\footnotesize $3L$}}

\put(0,12){\line(1,0){30}}
\put(98,42){\line(1,0){30}}

\thicklines
\put(0,30){\line(2,1){30}}
\put(20,45){\line(2,-1){30}}
\put(40,30){\line(2,1){30}}
\put(60,45){\line(2,-1){30}}
\put(80,30){\line(2,1){30}}

\end{picture}

\end{center}
\end{minipage}
\\
\hline

\phantom{$\surd$}\ $\textnormal{E}_7$ &
\begin{minipage}[m]{.25\linewidth}\setlength{\unitlength}{.20mm}

\begin{center}
\begin{picture}(500,60)(65,10)
\multiput(100,40)(30,0){6}{\circle*{5}} \put(160,20){\circle*{5}}
\put(100,40){\line(1,0){150}}
 \put(160,40){\line(0,-1){20}}
\put(72.5,40){\line(1,0){27.5}}
\put(68,38){\framebox(4,4)}

\put(100,45){$2$} \put(130,45){$3$} \put(160,45){$4$}
 \put(190,45){$3$}
\put(220,45){$2$}
 \put(150,20){$2$}
\end{picture}
\end{center}
\end{minipage}
&
\begin{minipage}[m]{.25\linewidth}\setlength{\unitlength}{.20mm}
\begin{center}
\begin{picture}(400,50)(-10,10)

\qbezier[15](45,42)(58,36)(71,30)
\put(41,45){\mbox{\footnotesize $2L$}}
\qbezier[15](12,50)(10,28)(10,15)
\put(16,20){\mbox{\footnotesize $L$}}

\put(117,33){\line(1,0){30}}

\thicklines
\put(0,30){\line(2,1){30}}
\put(20,45){\line(2,-1){30}}
\put(40,30){\line(2,1){30}}
\put(60,45){\line(2,-1){30}}
\put(80,30){\line(2,1){30}}
\put(100,45){\line(2,-1){30}}

\end{picture}
\end{center}
\end{minipage}
\\
\hline

\end{longtable}
\end{center}

\begin{center}
\begin{longtable}{|l|c|c|}
\hline \multicolumn{3}{|c|}
{Degree $3$} \\
\hline

\phantom{$\surd$\ $2\textnormal{A}_1+\textnormal{A}_2$}&Configuration& Example \\
\hline $\surd$\ $\textnormal{D}_5$&
\begin{minipage}[m]{.25\linewidth}\setlength{\unitlength}{.20mm}
\begin{center}
\begin{picture}(500,60)(90,5)
\multiput(130,40)(30,0){4}{\circle*{5}} \put(250,40){\circle{5}}
\put(160,20){\circle*{5}} \put(97.5,40){\circle{5}}
\put(100,40){\line(1,0){147.5}}
 \put(160,40){\line(0,-1){20}}
\put(100,45){$2$} \put(130,45){$3$} \put(160,45){$4$}
 \put(190,45){$3$}
\put(220,45){$2$}
\put(165,20){$2$}
\end{picture}
\end{center}
\end{minipage}

&
\begin{minipage}[m]{.25\linewidth}\setlength{\unitlength}{.20mm}
\begin{center}
\begin{picture}(400,50)(-10,10)

\qbezier[15](40,45)(27,33)(15,20)
\put(43,50){\mbox{\footnotesize $3L$}}

\put(78,32){\line(1,0){30}}
\put(0, 25){\line(1,0){30}} \thicklines
\put(0,30){\line(2,1){30}}
\put(20,45){\line(2,-1){30}}
\put(40,30){\line(2,1){30}}
\put(60,45){\line(2,-1){30}}

\end{picture}
\end{center}

\end{minipage}
\\

\hline $\surd$\  $\textnormal{A}_5$&
\begin{minipage}[m]{.25\linewidth}\setlength{\unitlength}{.20mm}

\begin{center}
\begin{picture}(500,60)(120,5)
\multiput(160,40)(30,0){4}{\circle*{5}} \put(160,20){\circle*{5}}
\put(127.5,40){\circle{5}} \put(130,40){\line(1,0){120}}
 \put(160,40){\line(0,-1){20}}
\put(130,45){$3$} \put(160,45){$4$}
 \put(190,45){$3$}
\put(220,45){$2$}
\put(165,20){$2$}
\end{picture}
\end{center}
\end{minipage}&
\begin{minipage}[m]{.25\linewidth}\setlength{\unitlength}{.20mm}
\begin{center}
\begin{picture}(400,50)(-10,10)

\qbezier[15](40,45)(30,37)(20,30)
\put(43,45){\mbox{\footnotesize $3L$}}

\put(98,42){\line(1,0){30}}

\thicklines
\put(0,30){\line(2,1){30}}
\put(20,45){\line(2,-1){30}}
\put(40,30){\line(2,1){30}}
\put(60,45){\line(2,-1){30}}
\put(80,30){\line(2,1){30}}

\end{picture}
\end{center}
\end{minipage}
\\

\hline \phantom{$\surd$}\ $\textnormal{E}_6$ &
\begin{minipage}[m]{.25\linewidth}\setlength{\unitlength}{.20mm}
\begin{center}
\begin{picture}(500,50)(90,10)
\multiput(130,40)(30,0){5}{\circle*{5}} \put(190,20){\circle*{5}}
\put(97.5,40){\circle{5}} \put(280,38){\framebox(4,4)}
\put(100,40){\line(1,0){180}}
 \put(190,40){\line(0,-1){20}}
\put(130,45){$2$}

\put(160,45){$3$} \put(190,45){$4$}
 \put(220,45){$3$}
\put(250,45){$2$}
\put(195,20){$2$}
\end{picture}
\end{center}

\end{minipage}&
\begin{minipage}[m]{.25\linewidth}\setlength{\unitlength}{.20mm}
\begin{center}
\begin{picture}(400,40)(-10,20)

\qbezier[15](45,42)(55,36)(65,30)
\put(41,45){\mbox{\footnotesize $2L$}}
\qbezier[15](5,50)(5,35)(5,20)
\put(10,20){\mbox{\footnotesize $L$}}

\put(98,42){\line(1,0){30}}

\thicklines
\put(0,30){\line(2,1){30}}
\put(20,45){\line(2,-1){30}}
\put(40,30){\line(2,1){30}}
\put(60,45){\line(2,-1){30}}
\put(80,30){\line(2,1){30}}

\end{picture}
\end{center}

\end{minipage}\\
\hline
\end{longtable}
\end{center}

\newpage

\begin{center}
\begin{longtable}{|l|c|c|}
\hline\multicolumn{3}{|c|}
{Degree $4$} \\
\hline

\phantom{$\surd$\ $2\textnormal{A}_1+\textnormal{A}_2$}&Configuration& Example \\

\hline $\surd$\  $\textnormal{D}_4$&
\begin{minipage}[m]{.25\linewidth}\setlength{\unitlength}{.20mm}
\begin{center}
\begin{picture}(500,60)(90,5)
\multiput(130,40)(30,0){3}{\circle*{5}} \put(220,40){\circle{5}}
\put(160,20){\circle*{5}} \put(100,40){\circle{5}}
\put(102.5,40){\line(1,0){115}}
 \put(160,40){\line(0,-1){20}}
\put(100,45){$2$} \put(130,45){$3$} \put(160,45){$4$}
 \put(190,45){$3$}
\put(220,45){$2$}  \put(165,20){$2$}
\end{picture}
\end{center}
\end{minipage}

&
\begin{minipage}[m]{.25\linewidth}\setlength{\unitlength}{.20mm}
\begin{center}
\begin{picture}(400,50)(-10,10)

\qbezier[15](40,45)(27,33)(15,20)
\put(45,45){\mbox{\footnotesize $3L$}}

\put(58,43){\line(1,0){30}}
\put(0, 25){\line(1,0){30}} \thicklines
\put(0,30){\line(2,1){30}}
\put(20,45){\line(2,-1){30}}
\put(40,30){\line(2,1){30}}

\end{picture}
\end{center}

\end{minipage}

\\

\hline $\surd$\ $\textnormal{A}_4$&
\begin{minipage}[m]{.25\linewidth}\setlength{\unitlength}{.20mm}
\begin{center}
\begin{picture}(500,60)(120,5)
\multiput(160,40)(30,0){3}{\circle*{5}} \put(250,40){\circle{5}}
\put(160,20){\circle*{5}} \put(127.5,40){\circle{5}}
\put(130,40){\line(1,0){117.5}}
 \put(160,40){\line(0,-1){20}}
 \put(130,45){$3$} \put(160,45){$4$}
 \put(190,45){$3$}
\put(220,45){$2$}
\put(165,20){$2$}
\end{picture}
\end{center}
\end{minipage}

&
\begin{minipage}[m]{.25\linewidth}\setlength{\unitlength}{.20mm}
\begin{center}
\begin{picture}(400,50)(-10,10)

\qbezier[30](0,48)(80,40)(15,20)
\put(50,35){\mbox{\footnotesize $3L$}}

\put(18,52){\line(1,0){30}}
\put(35, 23){\line(1,0){30}} \thicklines
\put(0,40){\line(2,1){30}}

\put(0,20){\line(2,1){30}}
\put(20,35){\line(2,-1){30}}

\end{picture}
\end{center}

\end{minipage}

\\
\hline $\surd$\ $\textnormal{A}_3+\textnormal{A}_1$&
\begin{minipage}[m]{.25\linewidth}\setlength{\unitlength}{.20mm}
\begin{center}
\begin{picture}(500,60)(150,10)
\multiput(220,40)(30,0){3}{\circle*{5}} \put(160,40){\circle*{5}}
\put(190,40){\circle{5}} \put(192.5,40){\line(1,0){90}}
 \put(160,40){\line(1,0){27.5}}
\put(190,45){$4$}
 \put(220,45){$3$}
\put(250,45){$2$}
\put(160,45){$2$}
\end{picture}
\end{center}
\end{minipage}&
\begin{minipage}[m]{.25\linewidth}\setlength{\unitlength}{.20mm}
\begin{center}
\begin{picture}(400,50)(-10,8)

\qbezier[15](5,45)(5,26)(5,12)
\put(4,45){\mbox{\footnotesize $L$}}
\qbezier[15](55,45)(55,30)(55,15)
\put(58,35){\mbox{\footnotesize $2L$}}

\put(37,32){\line(1,0){30}}
\put(0,20){\line(1,0){30}}
\put(0, 25){\line(1,0){30}} \thicklines
\put(0,30){\line(2,1){30}}
\put(20,45){\line(2,-1){30}}

\end{picture}
\end{center}
\end{minipage}
\\

\hline \phantom{$\surd$}\ $\textnormal{D}_5$&
\begin{minipage}[m]{.25\linewidth}\setlength{\unitlength}{.20mm}
\begin{center}
\begin{picture}(500,60)(90,5)
\multiput(130,40)(30,0){4}{\circle*{5}}
\put(247.5,37.5){\framebox(5,5)} \put(160,20){\circle*{5}}
\put(97.5,40){\circle{5}} \put(100,40){\line(1,0){147.5}}
 \put(160,40){\line(0,-1){20}}
\put(100,45){$2$} \put(130,45){$3$} \put(160,45){$4$}
 \put(190,45){$3$}
\put(220,45){$2$}
\put(165,20){$2$}
\end{picture}
\end{center}
\end{minipage}

&
\begin{minipage}[m]{.25\linewidth}\setlength{\unitlength}{.20mm}
\begin{center}
\begin{picture}(400,50)(-10,10)

\qbezier[15](5,45)(5,33)(5,15)
\put(10,20){\mbox{\footnotesize $L$}}
\qbezier[15](40,45)(55,33)(70,25)
\put(40,45){\mbox{\footnotesize $2L$}}

\put(78,32){\line(1,0){30}}

\thicklines
\put(0,30){\line(2,1){30}}
\put(20,45){\line(2,-1){30}}
\put(40,30){\line(2,1){30}}
\put(60,45){\line(2,-1){30}}

\end{picture}
\end{center}

\end{minipage}

\\
\hline

\end{longtable}
\end{center}

\begin{center}
\begin{longtable}{|l|c|c|}
\hline \multicolumn{3}{|c|}
{Degree $5$} \\
\hline

\phantom{$\surd$\ $2\textnormal{A}_1+\textnormal{A}_2$}&Configuration& Example \\

\hline $\surd$\ $\textnormal{A}_3$&
\begin{minipage}[m]{.25\linewidth}\setlength{\unitlength}{.20mm}
\begin{center}
\begin{picture}(500,60)(120,5)
\multiput(160,40)(30,0){2}{\circle*{5}} \put(220,40){\circle{5}}
\put(160,20){\circle*{5}} \put(127.5,40){\circle{5}}
\put(130,40){\line(1,0){87.5}}
 \put(160,40){\line(0,-1){20}}
 \put(130,45){$3$} \put(160,45){$4$}
 \put(190,45){$3$}
\put(220,45){$2$}
\put(165,20){$2$}
\end{picture}
\end{center}
\end{minipage}

&
\begin{minipage}[m]{.25\linewidth}\setlength{\unitlength}{.20mm}
\begin{center}
\begin{picture}(400,50)(-10,10)

\qbezier[30](0,48)(80,40)(15,20)
\put(20,50){\mbox{\footnotesize $3L$}}

\put(25,40){\line(1,0){30}}
\put(35, 23){\line(1,0){30}} \thicklines

\put(0,20){\line(2,1){30}}
\put(20,35){\line(2,-1){30}}

\end{picture}
\end{center}

\end{minipage}

\\
\hline $\surd$\ $\textnormal{A}_2+\textnormal{A}_1$&
\begin{minipage}[m]{.25\linewidth}\setlength{\unitlength}{.20mm}
\begin{center}
\begin{picture}(500,60)(150,10)
\multiput(220,40)(30,0){2}{\circle*{5}}
\put(160,40){\circle*{5}}\put(280,40){\circle{5}}
\put(190,40){\circle{5}} \put(192.5,40){\line(1,0){85}}
 \put(160,40){\line(1,0){27.5}}
\put(190,45){$4$}
 \put(220,45){$3$}
\put(250,45){$2$}
\put(160,45){$2$}
\end{picture}
\end{center}
\end{minipage}&
\begin{minipage}[m]{.25\linewidth}\setlength{\unitlength}{.20mm}
\begin{center}
\begin{picture}(400,50)(-10,8)

\qbezier[15](5,45)(5,26)(5,12)
\put(4,47){\mbox{\footnotesize $L$}}
\qbezier[15](55,45)(55,30)(55,15)
\put(58,45){\mbox{\footnotesize $2L$}}

\put(37,32){\line(1,0){30}}

\put(0, 20){\line(1,0){30}}

 \thicklines
\put(0,30){\line(2,1){30}}
\put(20,45){\line(2,-1){30}}

\end{picture}
\end{center}
\end{minipage}
\\

\hline \phantom{$\surd$}\  $\textnormal{A}_4$&\hspace{-13mm}
\begin{minipage}[m]{.25\linewidth}\setlength{\unitlength}{.20mm}
\begin{center}
\begin{picture}(500,60)(90,5)
\multiput(160,40)(30,0){3}{\circle*{5}}
\put(247.5,37.5){\framebox(4,4)}
\put(160,20){\circle*{5}}\put(130,40){\circle{5}}
 \put(132.5,40){\line(1,0){115}}
 \put(160,40){\line(0,-1){20}}
\put(130,45){$3$} \put(160,45){$4$}
 \put(190,45){$3$}
\put(220,45){$2$}
\put(165,20){$2$}
\end{picture}
\end{center}
\end{minipage}

&
\begin{minipage}[m]{.25\linewidth}\setlength{\unitlength}{.20mm}
\begin{center}
\begin{picture}(400,50)(-10,10)

\qbezier[15](5,45)(5,33)(5,15)
\put(12,12){\mbox{\footnotesize $L$}}
\qbezier[15](40,45)(55,33)(70,25)
\put(75,20){\mbox{\footnotesize $2L$}}

\put(56,42){\line(1,0){30}}

\thicklines
\put(0,30){\line(2,1){30}}
\put(20,45){\line(2,-1){30}}
\put(40,30){\line(2,1){30}}

\end{picture}
\end{center}

\end{minipage}

\\
\hline

\end{longtable}
\end{center}

\begin{center}
\begin{longtable}{|l|c|c|}
\hline \multicolumn{3}{|c|}
{Degree $6$} \\
\hline

\phantom{$\surd$\ $2\textnormal{A}_1+\textnormal{A}_2$}&Configuration& Example \\
\hline $\surd$\ $\textnormal{A}_2$&
\begin{minipage}[m]{.25\linewidth}\setlength{\unitlength}{.20mm}
\begin{center}
\begin{picture}(500,60)(120,5)
\multiput(160,40)(30,0){1}{\circle*{5}} \put(190,40){\circle{5}}
\put(160,20){\circle*{5}} \put(127.5,40){\circle{5}}
\put(130,40){\line(1,0){57.5}}
 \put(160,40){\line(0,-1){20}}
 \put(130,45){$3$} \put(160,45){$4$}
 \put(190,45){$3$}
\put(165,20){$2$}
\end{picture}
\end{center}
\end{minipage}

&
\begin{minipage}[m]{.25\linewidth}\setlength{\unitlength}{.20mm}
\begin{center}
\begin{picture}(400,50)(-10,10)

\qbezier[30](0,48)(80,40)(15,20)
\put(40,45){\mbox{\footnotesize $3L$}}

\put(35, 23){\line(1,0){30}} \thicklines

\put(0,20){\line(2,1){30}}
\put(20,35){\line(2,-1){30}}

\end{picture}
\end{center}

\end{minipage}

\\
\hline $\surd$\ $2\textnormal{A}_1$&
\begin{minipage}[m]{.25\linewidth}\setlength{\unitlength}{.20mm}
\begin{center}
\begin{picture}(500,60)(120,5)
\multiput(190,40)(30,0){1}{\circle*{5}} \put(220,40){\circle{5}}
 \put(127.5,40){\circle*{5}}\put(160,40){\circle{5}}
\put(130,40){\line(1,0){27.5}}\put(162.5,40){\line(1,0){55}}
 \put(130,45){$2$} \put(160,45){$4$}
 \put(190,45){$3$}
\put(220,45){$2$}
\end{picture}
\end{center}
\end{minipage}

&
\begin{minipage}[m]{.25\linewidth}\setlength{\unitlength}{.20mm}
\begin{center}
\begin{picture}(400,50)(-10,10)

\qbezier[15](40,50)(40,35)(40,20)
\put(20,45){\mbox{\footnotesize $3L$}}

\put(25,40){\line(1,0){30}}
\put(18, 33){\line(1,0){40}} \thicklines

\put(0,20){\line(2,1){30}}

\end{picture}
\end{center}

\end{minipage}

\\

\hline  \phantom{$\surd$}\ $\textnormal{A}_2+\textnormal{A}_1$&
\begin{minipage}[m]{.25\linewidth}\setlength{\unitlength}{.20mm}
\begin{center}
\begin{picture}(500,60)(90,5)
\multiput(130,40)(30,0){1}{\circle*{5}}
\multiput(190,40)(30,0){1}{\circle{5}}
\put(160,40){\circle*{5}}\put(220,40){\circle*{5}}
\put(102.5,40){\line(1,0){55}} \put(162.5,40){\line(1,0){25}}
 \put(192.5,40){\line(1,0){25}}
\put(130,45){$2$} \put(160,45){$3$}
\put(190,45){$4$}\put(220,45){$2$}
\put(98,38){\framebox(4,4)}
\end{picture}
\end{center}

\end{minipage}

&
\begin{minipage}[m]{.25\linewidth}\setlength{\unitlength}{.20mm}
\begin{center}
\begin{picture}(400,50)(-10,10)

\qbezier[15](5,45)(5,33)(5,15)
\put(10,15){\mbox{\footnotesize $L$}}
\qbezier[15](40,45)(55,33)(70,25)
\put(45,45){\mbox{\footnotesize $2L$}}

\put(37,34){\line(1,0){30}}

\thicklines
\put(0,30){\line(2,1){30}}
\put(20,45){\line(2,-1){30}}

\end{picture}
\end{center}

\end{minipage}

\\

\hline
\end{longtable}
\end{center}

\begin{center}
\begin{longtable}{|l|c|c|}
\hline \multicolumn{3}{|c|}
{Degree $7$} \\
\hline

\phantom{$\surd$\ $2\textnormal{E}_1+\textnormal{E}_1$}&Configuration& Example \\
\hline  $\surd$\ $\textnormal{A}_1$&
\begin{minipage}[m]{.25\linewidth}\setlength{\unitlength}{.20mm}
\begin{center}
\begin{picture}(500,60)(90,5)
\multiput(100,40)(30,0){1}{\circle*{5}}\multiput(130,40)(30,0){2}{\circle{5}}
 \put(100,40){\line(1,0){27.5}}\put(132.5,40){\line(1,0){25}}
 \put(100,45){$2$} \put(130,45){$4$}
\put(160,45){$3$}
\end{picture}
\end{center}
\end{minipage}

&
\begin{minipage}[m]{.25\linewidth}\setlength{\unitlength}{.20mm}
\begin{center}
\begin{picture}(400,50)(-10,10)

\qbezier[15](35,48)(45,38)(55,28)
\put(60,20){\mbox{\footnotesize $3L$}}

\put(18,42){\line(1,0){35}}
\thicklines
\put(0,30){\line(2,1){30}}

\end{picture}
\end{center}

\end{minipage}

\\

\hline
\end{longtable}
\end{center}
\normalsize
\end{Prop}
\pf Starting from the $\textnormal{E}_7$ dual graph for the degree $1$ case,  we apply
$\textbf{P-1}$ and $\textbf{P-2}$
successively to get the possible dual graphs. Then from the obtained possible dual graphs, we exclude the dual graphs that violate the properties in Lemma~\ref{wt4}.
\qed

\begin{Lem}\label{wt3}
 If the dual graph for an effective anticanonical divisor  on a weak del Pezzo surface  has a vertex of weight $3$ as a maximum, the divisor $\tilde{D}$ does not have two $0$-curves. In addition, if the
divisor has a $1$-curve, then its dual graph is  obtained  by suitable blow ups from the following:
\begin{center}
\begin{minipage}[m]{.2\linewidth}\setlength{\unitlength}{.28mm}\setlength{\unitlength}{.25mm}
\begin{center}
\begin{picture}(50,60)(100,5)

\multiput(130,40)(30,0){1}{\circle*{5}}
\multiput(190,40)(30,0){1}{\circle*{5}}
\put(160,40){\circle{5}}\put(220,40){\circle*{5}}
\put(102.5,40){\line(1,0){55}} \put(162.5,40){\line(1,0){55}}
\put(130,45){$2$} \put(160,45){$3$} \put(190,45){$2$}
\put(98,38){\framebox(4,4)}\put(100,40){\circle*{4}}
\end{picture}
\end{center}
\end{minipage}
\end{center}

\end{Lem}
\pf The proof is similar as that of Lemma~\ref{wt4}.  In this case, the possible sequences for~$\{\mult_{p_{i, k}}(C)\}_{k=1}^{j}$ have the form as follows:

$$(2,2,2, *, *,\cdots), \ (3,2, *, *,\cdots).$$
Therefore, the curve $C$ consists of only lines.
We also see that the first form must be $(2,2,2)$ and the second must be $(3,2)$.

For the case $(2,2,2)$, the curve $C$ must consists of one double line and one single line. We also see that one blow up at a point on the double line but not on the single line must be taken. Then one blow up at the intersection point of the exceptional divisor and the strict transform of the double line must follow. And then one more blow up at the intersection point of the exceptional divisor  of the second blow up and the strict transform of the double line must be taken.

For the case $(3,2)$, the curve $C$ must also consists of one double line and one single line.  We must take  the blow up at the intersection point of the double line and the single line. Then the blow up at a point on the strict transform of the double line but not on that of the single line must follow.

In both cases, the dual graph does not contain two $0$-curves. If it contains a $1$-curve, such a dual graph must be obtained in the way of the case $(2,2,2)$.
\qed

\begin{Prop}\label{proposition:wt3}
If the dual graph for an effective anticanonical divisor  on a weak del Pezzo surface  has a vertex of weight $3$ as a maximum,
then it is exactly one of the following:

\footnotesize
\begin{center}
\begin{longtable}{|l|c|c|}
\hline \multicolumn{3}{|c|}
{Degree $1$} \\
\hline

\phantom{$\surd$\ $2\textnormal{E}_1+\textnormal{E}_1$}&Configuration& Example \\
\hline $\surd$\ $\textnormal{E}_6$
&\begin{minipage}[m]{.25\linewidth}\setlength{\unitlength}{.20mm}

\begin{center}
\begin{picture}(500,70)(95,-10)
\multiput(100,40)(30,0){5}{\circle*{5}} \put(160,20){\circle*{5}}
\put(100,40){\line(1,0){120}} \put(160,40){\line(0,-1){37.5}}
\put(160,0){\circle{5}}
\put(130,45){$2$} \put(160,45){$3$} \put(190,45){$2$}
\put(150,20){$2$} 
\end{picture}
\end{center}
\end{minipage}
&
\begin{minipage}[m]{.25\linewidth}\setlength{\unitlength}{.20mm}
\begin{center}
\begin{picture}(400,70)(-10,-10)

\qbezier[30](30,0)(10,18)(100,20)
\put(33,-5){\mbox{\footnotesize $L$}}
\qbezier[30](40,45)(70,22)(100,0)
\put(41,46){\mbox{\footnotesize $2L$}}

\put(18,12){\line(1,0){30}}
\put(98,42){\line(1,0){30}}

\thicklines
\put(0,30){\line(2,1){30}}
\put(20,45){\line(2,-1){30}}
\put(40,30){\line(2,1){30}}
\put(60,45){\line(2,-1){30}}
\put(80,30){\line(2,1){30}}

\put(0,0){\line(2,1){30}}

\end{picture}
\end{center}
\end{minipage}
\\

\hline
\end{longtable}
\end{center}

\begin{center}
\begin{longtable}{|l|c|c|}
\hline \multicolumn{3}{|c|}
{Degree $2$} \\
\hline

\phantom{$\surd$\ $2\textnormal{E}_1+\textnormal{E}_1$}&Configuration& Example \\
\hline  $\surd$\ $\textnormal{D}_5$ &\hspace{-13mm}
\begin{minipage}[m]{.25\linewidth}\setlength{\unitlength}{.20mm}
\begin{center}
\begin{picture}(500,70)(60,10)
\multiput(130,60)(30,0){4}{\circle*{5}} \put(160,40){\circle*{5}}
\put(97.5,60){\circle{5}} \put(160,17.5){\circle{5}}
\put(100,60){\line(1,0){120}} \put(160,60){\line(0,-1){40}}
\put(130,65){$2$} \put(160,65){$3$}
\put(190,65){$2$} 
\put(150,38){$2$}
\end{picture}
\end{center}

\end{minipage}&
\begin{minipage}[m]{.25\linewidth}\setlength{\unitlength}{.20mm}
\begin{center}
\begin{picture}(400,50)(-10,10)

\qbezier[15](5,40)(10,25)(15,10)
\put(0,42){\mbox{\footnotesize $L$}}
\qbezier[15](25,60)(45,40)(26,30)
\put(40,40){\mbox{\footnotesize $2L$}}

\put(77,33){\line(1,0){30}}
\put(10,55){\line(1,0){30}}
\put(0,17){\line(1,0){30}}

\thicklines
\put(0,30){\line(2,1){30}}
\put(20,45){\line(2,-1){30}}
\put(40,30){\line(2,1){30}}
\put(60,45){\line(2,-1){30}}

\end{picture}
\end{center}
\end{minipage}\\
\hline
\end{longtable}
\end{center}

\newpage

\begin{center}
\begin{longtable}{|l|c|c|}
\hline \multicolumn{3}{|c|}
{Degree $3$} \\
\hline

\phantom{$\surd$\ $2\textnormal{E}_1+\textnormal{E}_1$}&Configuration& Example \\
\hline  $\surd$\ $\textnormal{D}_4$ &\hspace{-13mm}
\begin{minipage}[m]{.25\linewidth}\setlength{\unitlength}{.20mm}
\begin{center}
\begin{picture}(500,70)(60,10)
\multiput(130,60)(30,0){3}{\circle*{5}} \put(220,60){\circle{5}}
\put(160,40){\circle*{5}} \put(97.5,60){\circle{5}}
\put(160,17.5){\circle{5}} \put(100,60){\line(1,0){117.5}}
\put(160,60){\line(0,-1){40}}
\put(130,65){$2$} \put(160,65){$3$} \put(190,65){$2$}
\put(150,38){$2$}
\end{picture}
\end{center}

\end{minipage}&
\begin{minipage}[m]{.25\linewidth}\setlength{\unitlength}{.20mm}
\begin{center}
\begin{picture}(400,50)(-10,10)

\qbezier[15](5,40)(10,25)(15,10)
\put(5,42){\mbox{\footnotesize $L$}}
\qbezier[15](25,60)(45,40)(26,30)
\put(40,40){\mbox{\footnotesize $2L$}}

\put(57,43){\line(1,0){30}}
\put(10,55){\line(1,0){30}}
\put(0,17){\line(1,0){30}}

\thicklines
\put(0,30){\line(2,1){30}}
\put(20,45){\line(2,-1){30}}
\put(40,30){\line(2,1){30}}

\end{picture}
\end{center}
\end{minipage}\\
\hline
\end{longtable}
\end{center}

\begin{center}
\begin{longtable}{|l|c|c|}
\hline \multicolumn{3}{|c|}
{Degree $4$} \\
\hline

\phantom{$\surd$\ $2\textnormal{E}_1+\textnormal{E}_1$}&Configuration& Example \\

\hline  $\surd$\  $\textnormal{A}_3$&
\begin{minipage}[m]{.25\linewidth}\setlength{\unitlength}{.20mm}
\begin{center}
\begin{picture}(500,60)(90,5)
\multiput(130,40)(30,0){3}{\circle*{5}} \put(130,17.5){\circle{5}}
\put(102.5,40){\line(1,0){90}} \put(130,40){\line(0,-1){20}}
\put(100,45){$2$} \put(130,45){$3$} \put(160,45){$2$}
\put(135,20){$2$} \put(100,40){\circle{5}}
\end{picture}
\end{center}
\end{minipage}

&
\begin{minipage}[m]{.25\linewidth}\setlength{\unitlength}{.20mm}
\begin{center}
\begin{picture}(400,50)(-10,10)

\qbezier[15](0,40)(25,25)(0,10)
\put(0,42){\mbox{\footnotesize $3L$}}

\put(0, 20){\line(1,0){30}}
\put(38,33){\line(1,0){30}}
\put(0, 25){\line(1,0){30}} \thicklines
\put(0,30){\line(2,1){30}}
\put(20,45){\line(2,-1){30}}

\end{picture}
\end{center}

\end{minipage}

\\
\hline  $\surd$\  $\textnormal{A}_3$ &\hspace{-13mm}
\begin{minipage}[m]{.25\linewidth}\setlength{\unitlength}{.20mm}
\begin{center}
\begin{picture}(500,70)(60,10)
\multiput(130,60)(30,0){2}{\circle*{5}} \put(190,60){\circle{5}}
\put(160,40){\circle*{5}} \put(97.5,60){\circle{5}}
\put(160,17.5){\circle{5}} \put(100,60){\line(1,0){87.5}}
\put(160,60){\line(0,-1){40}}
\put(130,65){$2$} \put(160,65){$3$} \put(190,65){$2$}
\put(150,38){$2$}
\end{picture}
\end{center}

\end{minipage}&
\begin{minipage}[m]{.25\linewidth}\setlength{\unitlength}{.20mm}
\begin{center}
\begin{picture}(400,50)(-10,10)

\qbezier[10](5,40)(10,25)(15,10)
\put(5,42){\mbox{\footnotesize $L$}}
\qbezier[15](25,60)(50,40)(26,30)
\put(33,55){\mbox{\footnotesize $2L$}}

\put(57,43){\line(1,0){30}}

\put(0,17){\line(1,0){30}}

\thicklines
\put(0,30){\line(2,1){30}}
\put(20,45){\line(2,-1){30}}
\put(40,30){\line(2,1){30}}

\end{picture}
\end{center}
\end{minipage}\\
\hline
\end{longtable}
\end{center}

\begin{center}
\begin{longtable}{|l|c|c|}
\hline \multicolumn{3}{|c|}
{Degree $5$} \\
\hline

\phantom{$\surd$\ $2\textnormal{E}_1+\textnormal{E}_1$}&Configuration& Example \\

\hline  $\surd$\ $\textnormal{A}_2$&\hspace{-13mm}
\begin{minipage}[m]{.25\linewidth}\setlength{\unitlength}{.20mm}
\begin{center}
\begin{picture}(500,60)(90,5)
\multiput(130,40)(30,0){1}{\circle{5}}
\multiput(190,40)(30,0){2}{\circle*{5}} \put(160,40){\circle{5}}
 \put(162.5,40){\line(1,0){60}}\put(132.5,40){\line(1,0){25}}
\put(130,45){$2$} \put(160,45){$3$} \put(190,45){$2$}
\
\end{picture}
\end{center}
\end{minipage}

&
\begin{minipage}[m]{.25\linewidth}\setlength{\unitlength}{.20mm}
\begin{center}
\begin{picture}(400,50)(-10,10)
\qbezier[15](4,50)(10,35)(15,20)
\put(7,45){\mbox{\footnotesize $3L$}}

\put(0,39){\line(1,0){30}}
\put(30, 25){\line(1,0){30}} \thicklines

\put(0,20){\line(2,1){30}}
\put(20,35){\line(2,-1){30}}

\end{picture}
\end{center}

\end{minipage}\\

\hline  $\surd$\ $\textnormal{A}_2$&
\begin{minipage}[m]{.25\linewidth}\setlength{\unitlength}{.20mm}
\begin{center}
\begin{picture}(500,60)(90,5)
\multiput(130,40)(30,0){2}{\circle*{5}} \put(130,17.5){\circle{5}}
\put(102.5,40){\line(1,0){85}} \put(130,40){\line(0,-1){20}}
\put(100,45){$2$} \put(130,45){$3$} \put(160,45){$2$}
\put(135,20){$2$} \put(100,40){\circle{5}}\put(190,40){\circle{5}}
\end{picture}
\end{center}
\end{minipage}

&
\begin{minipage}[m]{.25\linewidth}\setlength{\unitlength}{.20mm}
\begin{center}
\begin{picture}(400,50)(-10,10)

\qbezier[15](5,45)(15,27)(5,10)
\put(0,49){\mbox{\footnotesize $3L$}}

\put(0, 20){\line(1,0){30}}
\put(18,42){\line(1,0){30}}
\put(0, 25){\line(1,0){30}} \thicklines
\put(0,30){\line(2,1){30}}

\end{picture}
\end{center}

\end{minipage}

\\
\hline
\end{longtable}
\end{center}

\begin{center}
\begin{longtable}{|l|c|c|}
\hline \multicolumn{3}{|c|}
{Degree $6$} \\
\hline

\phantom{$\surd$\ $2\textnormal{E}_1+\textnormal{E}_1$}&Configuration& Example \\

\hline  $\surd$\ $\textnormal{A}_1$&
\begin{minipage}[m]{.25\linewidth}\setlength{\unitlength}{.20mm}
\begin{center}
\begin{picture}(500,60)(90,5)
\multiput(160,40)(30,0){1}{\circle*{5}}\multiput(190,40)(30,0){1}{\circle{5}}
\put(102.5,40){\line(1,0){25}}\put(132.5,40){\line(1,0){55}}
\put(100,45){$2$} \put(130,45){$3$} \put(160,45){$2$}
 \put(100,40){\circle{5}}\put(130,40){\circle{5}}
\end{picture}
\end{center}
\end{minipage}

&
\begin{minipage}[m]{.25\linewidth}\setlength{\unitlength}{.20mm}
\begin{center}
\begin{picture}(400,50)(-10,10)
\qbezier[15](4,45)(10,25)(15,20)
\put(8,45){\mbox{\footnotesize $3L$}}

\put(0,39){\line(1,0){30}}
\put(18, 33){\line(1,0){30}} \thicklines

\put(0,20){\line(2,1){30}}

\end{picture}
\end{center}

\end{minipage}\\

\hline  $\surd$\ $\textnormal{A}_1$&
\begin{minipage}[m]{.25\linewidth}\setlength{\unitlength}{.20mm}
\begin{center}
\begin{picture}(500,60)(90,5)
\multiput(130,40)(30,0){1}{\circle*{5}} \put(130,17.5){\circle{5}}
\put(102.5,40){\line(1,0){55}} \put(130,40){\line(0,-1){20}}
\put(100,45){$2$} \put(130,45){$3$} \put(160,45){$2$}
\put(135,20){$2$} \put(100,40){\circle{5}}\put(160,40){\circle{5}}
\end{picture}
\end{center}
\end{minipage}

&
\begin{minipage}[m]{.25\linewidth}\setlength{\unitlength}{.20mm}
\begin{center}
\begin{picture}(400,50)(-10,10)
\qbezier[15](5,39)(30,27)(5,15)
\put(3,42){\mbox{\footnotesize $3L$}}

\put(0, 20){\line(1,0){30}}
\put(0,30){\line(1,0){30}}
\put(0, 25){\line(1,0){30}} \thicklines

\end{picture}
\end{center}

\end{minipage}

\\
\hline
\end{longtable}
\end{center}

\begin{center}
\begin{longtable}{|l|c|c|}
\hline \multicolumn{3}{|c|}
{Degree $7$} \\
\hline

\phantom{$\surd$\ $2\textnormal{E}_1+\textnormal{E}_1$}&Configuration& Example \\
\hline \phantom{$\surd$}\  $\textnormal{A}_1$ &
\begin{minipage}[m]{.25\linewidth}\setlength{\unitlength}{.20mm}
\begin{center}
\begin{picture}(500,60)(90,10)

\put(102.5,40){\line(1,0){25} }
\put(132.5,40){\line(1,0){25}} \put(162.5,40){\line(1,0){25}}
\multiput(100,40)(90,0){1}{\circle{5}}\multiput(160,40)(90,0){1}{\circle*{5}}
\multiput(130,40)(30,0){1}{\circle{5}}\put(188,38){\framebox(4,4)}
\multiput(160,45)(30,0){2}{$2$}\multiput(100,45)(30,0){1}{$2$}\multiput(130,45)(30,0){1}{$3$}
\end{picture}
\end{center}

\end{minipage}&
\begin{minipage}[m]{.25\linewidth}\setlength{\unitlength}{.20mm}
\begin{center}
\begin{picture}(400,50)(-10,15)

\qbezier[12](5,40)(10,30)(15,20)
\put(0,42){\mbox{\footnotesize $L$}}
\qbezier[15](25,60)(45,40)(26,30)
\put(40,50){$2L$}

\put(17,42){\line(1,0){30}}

\thicklines
\put(0,30){\line(2,1){30}}

\end{picture}
\end{center}
\end{minipage}\\
\hline \phantom{$\surd$}\ $\textnormal{A}_1$&
\begin{minipage}[m]{.25\linewidth}\setlength{\unitlength}{.20mm}
\begin{center}
\begin{picture}(500,55)(120,12)
 \put(190,40){\circle{5}} \put(128,38){\framebox(4,4)}
 \put(160,40){\circle*{5}}
\put(132.5,40){\line(1,0){25}} \put(162,40){\line(1,0){25}}
\put(159,45){$2$}\put(129,45){$3$} \put(189,45){$2$}
\end{picture}
\end{center}
\end{minipage}

&
\begin{minipage}[m]{.25\linewidth}\setlength{\unitlength}{.20mm}
\begin{center}
\begin{picture}(400,50)(-10,10)

\qbezier[15](4,40)(10,25)(15,15)
\put(4,45){\mbox{\footnotesize $3L$}}

\put(16, 33){\line(1,0){30}} \thicklines

\put(0,20){\line(2,1){30}}

\end{picture}
\end{center}

\end{minipage}\\
\hline
\end{longtable}
\end{center}
\normalsize
\end{Prop}
\pf Starting from the $\textnormal{E}_6$ dual graph for the degree $1$ case,  we apply
$\textbf{P-1}$ and $\textbf{P-2}$
successively to get the possible dual graphs. Then from the obtained possible dual graphs, we exclude the dual graphs that violate the properties in Lemma~\ref{wt3}.\qed

For the dual graphs of effective anticanonical divisors on  weak
del Pezzo surfaces  with a vertex of weight $2$ as a maximum, let
$\tilde{D}$ be such an anticanonical divisor. As in
formula~(\ref{adjuction}), the divisor $\tilde{D}$ is of the form
\begin{equation*}\sum a_h\tilde{C_h}+\sum_{i, j}\left(\sum_{k=1}^j\left(\mult_{p_{i,k}}(C)-1\right)\right)E_{i,j}\end{equation*}
with the same notation for formula~(\ref{adjuction}).

For some $(i, j)$, we must have $
\sum_{k=1}^j\left(\mult_{p_{i,k}}(C)-1\right)=2 $. The possible sequences for~$\{\mult_{p_{i, k}}(C)\}_{k=1}^{j}$ have the form as follows:

$$(2,2, *, *,\cdots), \ (3, *, *,\cdots).$$
Furthermore, we can see that only the sequences
\[(2,2,1,1,1,1), \ (2,2,1,1,1), \  (2,2,1,1), \ (2,2,1)\]
\[  (2,2),\]
\[(3,1,1), \ (3,1), \ (3)\]
can happen.
For the four sequences in the first row, the curve $C$ must consist of one irreducible conic and a line intersecting tangentially. For the  sequence in the second row, the curve $C$ consists either of one irreducible conic and a line intersecting tangentially or of one double line and one single line. For the three sequences in the last row, the curve $C$ consists either of one double line and one single line or of three lines intersecting at a single point.

 When the curve $C$ is given with one of the sequences above, the way to take blow ups is unique except the case
where $C$ consists of one irreducible conic and a line intersecting tangentially and the sequence is $(2,2,1)$. This exceptional case has two ways to take blow ups. First, we take the blow up at the intersection point of the conic and the line. Then the blow up at the intersection point of the strict transforms of the conic and the line follows. For the last blow up, we have two choices. One is to take the blow up at the intersection point of the exceptional divisor of the second blow up and the strict transform of the conic, and the other is  to take the blow up at the intersection point of the exceptional divisor of the second blow up and the strict transform of the line.

\begin{Lem}\label{wt2-1}
 If the divisor $\tilde{D}$ has either a $1$-curve or a $2$-curve, it cannot have any other
curve with nonnegative self-intersection number.
\end{Lem}
\pf For the divisor $\tilde{D}$ to have either a $1$-curve or a $2$-curve, the sequence must be either~$(2,2,1)$ or $(2, 2)$.
The curve $C$ consists either of one irreducible conic and a line intersecting tangentially or of one double line and one single line. If the curve $C$ consists  of one irreducible conic and a line intersecting tangentially, then the conic becomes either a $1$-curve or a $2$-curve and the single line becomes either a $-1$-curve or a $-2$-curve. If the curve $C$ consists of one double line and one single line, then the single line becomes a $1$-curve and the double line becomes a $-1$-curve.
 \qed

\begin{Lem}\label{wt2-2}
If the dual graph of $\widetilde{D}$ contains a chain
consisting of five or four vertices with weight~$2$ and negative self-intersection, then it is  obtained  by suitable blow ups from the following:
\begin{center}
\begin{minipage}[m]{.2\linewidth}\setlength{\unitlength}{.28mm}\setlength{\unitlength}{.25mm}
\begin{center}
\begin{picture}(100,60)(200,5)
\put(100,40){\line(1,0){27.5} }
\put(132.5,40){\line(1,0){25}} \put(162.5,40){\line(1,0){25}}
\put(247.5,40){\line(-1,0){25}}\put(190,40){\line(1,0){27.5}}
\multiput(100,40)(150,0){1}{\circle*{5}}
\multiput(160,40)(30,0){2}{\circle*{5}}
\multiput(130,40)(90,0){2}{\circle{5}}\multiput(250,40)(90,0){1}{\circle{5}}
\multiput(130,45)(30,0){4}{$2$}
\end{picture}
\end{center}
\end{minipage}
\begin{minipage}[m]{.2\linewidth}\setlength{\unitlength}{.28mm}\setlength{\unitlength}{.25mm}
\begin{center}
\begin{picture}(100,60)(100,5)
\put(100,40){\line(1,0){27.5} }
\put(132.5,40){\line(1,0){25}} \put(162.5,40){\line(1,0){25}}
\put(280,40){\line(-1,0){27.5}}\put(190,40){\line(1,0){27.5}}\put(220,40){\line(1,0){27.5}}
\multiput(100,40)(180,0){2}{\circle*{5}}
\multiput(160,40)(30,0){3}{\circle*{5}}
\multiput(130,40)(120,0){2}{\circle{5}}
\multiput(130,45)(30,0){5}{$2$}
\end{picture}
\end{center}
\end{minipage}
\end{center}

\end{Lem}
\pf If the dual graph of the divisor $\tilde{D}$ contains a chain
consisting of five or four vertices with weight $2$ and negative self-intersection, then we have two possibilities.
 One possibility is that we have either $(2,2,1,1,1,1)$ or  $(2,2,1,1,1)$ for the sequence $\{\mult_{p_{i, k}}(C)\}_{k=1}^{j}$ for some $(i,j)$. In this case, the assertion is clear.

The other possibility is as follows: the curve $C$  consists of one double line and one single line; for  the sequence $\{\mult_{p_{i, k}}(C)\}_{k=1}^{j}$, $(3,1,1)$, $(3,1)$, or $(3)$ is attained over the intersection point of the double line and the single line; the sequence $(2,2)$ is attained over a point on the double line but not on the single line.  Then, the blow ups in the way corresponding to each sequence complete the proof.
 \qed

\begin{Lem}\label{wt2-3} For each of $k=0, 1$, consider the set of dual graphs of all effective anticanonical divisors on weak del Pezzo surfaces such that they
have exactly one vertex with weight $1$ and  self-intersection number $k$  as a
maximal self-intersection number.
If the dual graph of the divisor~$\tilde{D}$ has exactly one vertex with weight $1$ and  self-intersection number $k$  as a
maximum self-intersection number and has  a longest
chain consisting of vertices of weight 2 in the set, then it is  obtained  by suitable blow ups from the following:
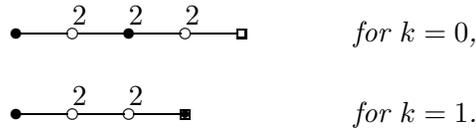
\begin{center}
\begin{minipage}[m]{.2\linewidth}\setlength{\unitlength}{.28mm}\setlength{\unitlength}{.25mm}
\begin{center}
\begin{picture}(100,60)(160,10)

\put(100,40){\line(1,0){27.5} }
\put(132.5,40){\line(1,0){25}} \put(162.5,40){\line(1,0){25}}
\multiput(100,40)(60,0){2}{\circle*{5}}\put(192.5,40){\line(1,0){25}}
\multiput(130,40)(60,0){1}{\circle{5}}\put(218,38){\framebox(4,4)}\put(190,40){\circle{5}}
\multiput(130,45)(30,0){3}{$2$}
\end{picture}
\end{center}
\end{minipage}
for $k=0$,
\end{center}

\vspace{-12mm}
\begin{center}
\begin{minipage}[m]{.2\linewidth}\setlength{\unitlength}{.28mm}\setlength{\unitlength}{.25mm}
\begin{center}
\begin{picture}(100,60)(160,10)
\put(100,40){\line(1,0){27.5} }
\put(132.5,40){\line(1,0){25}} \put(162.5,40){\line(1,0){25}}
\multiput(100,40)(90,0){1}{\circle*{5}}
\multiput(130,40)(30,0){2}{\circle{5}}\put(188,38){\framebox(4,4)}\put(190,40){\circle*{4}}
\multiput(130,45)(30,0){2}{$2$}
\end{picture}
\end{center}
\end{minipage}
for $k=1$.
\end{center}

\end{Lem}
\pf For $k=0$, only in the following cases the dual graph of $\tilde{D}$ can satisfy the required conditions.
The first case is when we have the sequence $(2,2,1,1)$. The other case is when the curve $C$  consists of one double line and one single line; for  the sequence $\{\mult_{p_{i, k}}(C)\}_{k=1}^{j}$, $(3)$ is attained over the intersection point of the double line and the single line; the sequence $(2,2)$ is attained over a point on the double line but not on the single line.
The blow ups in the way corresponding to each sequence complete the proof for the case $k=0$.

 For $k=1$, only in  the following cases the dual graph of $\tilde{D}$ can satisfy the required conditions.
The first case is when we have the sequence $(2,2,1)$. The other case is when the curve $C$ consists of one double line and one single line with the sequence $(2,2)$.
 \qed
\begin{Lem}\label{wt2-4}
If the divisor $\tilde{D}$ has two reduced  $0$-curves, three
reduced $0$-curves, or exactly one $2$-curve as a maximal self-intersection
number, then its dual graph is   obtained  by suitable blow ups
from the following:
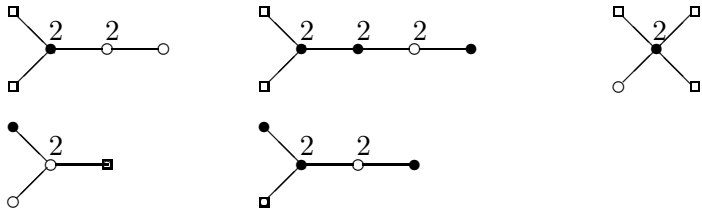
\newline

\begin{minipage}[m]{.2\linewidth}\setlength{\unitlength}{.28mm}\setlength{\unitlength}{.25mm}
\begin{center}
\begin{picture}(50,60)(125,5)
\put(132,58){\line(1,-1){20}} \put(132,22){\line(1,1){18}}
\put(150,40){\line(1,0){27.5} }
\put(182.5,40){\line(1,0){25} }
\multiput(150,40)(30,0){1}{\circle*{5}}
\put(128,58){\framebox(4,4)}
\multiput(128,18)(100,0){1}{\framebox(4,4)}
\put(180,40){\circle{5}} \put(210,40){\circle{5}}
\multiput(149,45)(30,0){2}{$2$}
\end{picture}
\end{center}
\end{minipage}
\begin{minipage}[m]{.2\linewidth}\setlength{\unitlength}{.28mm}\setlength{\unitlength}{.25mm}
\begin{center}
\begin{picture}(50,60)(125,5)
\put(132,58){\line(1,-1){18}} \put(132,22){\line(1,1){20}}
\put(150,40){\line(1,0){57.5} }
\put(212.5,40){\line(1,0){25}} \put(240,40){\circle*{5}}
\multiput(150,40)(30,0){2}{\circle*{5}} 
\multiput(128,18)(100,0){1}{\framebox(4,4)}\put(210,40){\circle{5}}\multiput(128,58)(100,0){1}{\framebox(4,4)}
\multiput(149,45)(30,0){3}{$2$}
\end{picture}
\end{center}
\end{minipage}
\begin{minipage}[m]{.2\linewidth}\setlength{\unitlength}{.28mm}\setlength{\unitlength}{.25mm}
\begin{center}
\begin{picture}(50,60)(70,5)
\multiput(128,58)(100,0){1}{\framebox(4,4)}
\multiput(130,20)(100,0){1}{\circle{5}}
 \put(150,40){\circle*{5}}
\put(168.5,58){\framebox(4,4)} \put(168.5,18){\framebox(4,4)}
\put(132,58){\line(1,-1){20}} \put(132.5,22.5){\line(1,1){20}}
\put(150,40){\line(1,1){18.3}} \put(150,40){\line(1,-1){18.3}}
\put(148,45){$2$}
\end{picture}
\end{center}
\end{minipage}

\begin{minipage}[m]{.2\linewidth}\setlength{\unitlength}{.28mm}\setlength{\unitlength}{.25mm}
\begin{center}
\begin{picture}(50,60)(125,5)
\put(130,60){\line(1,-1){18}} \put(132,22){\line(1,1){16.2}}

\put(152.5,40){\line(1,0){28} }

\put(178,38){\framebox(4,4)} \put(180,40){\circle{4}}

\multiput(150,40)(30,0){1}{\circle{5}}
\put(130,60){\circle*{5}}
\multiput(130,20)(100,0){1}{\circle{5}}

\multiput(149,45)(30,0){1}{$2$}
\end{picture}
\end{center}
\end{minipage}
\begin{minipage}[m]{.2\linewidth}\setlength{\unitlength}{.28mm}\setlength{\unitlength}{.25mm}
\begin{center}
\begin{picture}(50,60)(125,5)
\put(132,58){\line(1,-1){20}} \put(132,22){\line(1,1){18}}
\put(150,40){\line(1,0){27.5} }
\put(182.5,40){\line(1,0){25} }
\multiput(150,40)(30,0){1}{\circle*{5}}
 \put(180,40){\circle{5}}
\put(210,40){\circle*{5}}\put(128,18){\framebox(4,4)}\put(130,60){\circle*{5}}\put(130,20){\circle{4}}
\multiput(149,45)(30,0){2}{$2$}
\end{picture}
\end{center}
\end{minipage}

\end{Lem}
\pf For the divisor $\tilde{D}$ to have either two reduced $0$-curves or
three reduced $0$-curves, the curve $C$ consists  of three lines intersecting at a
single point with  the sequence $(3,1,1)$, $(3, 1)$ or
$(3)$. If the divisor $\tilde{D}$ has exactly one $2$-curve, then
the curve $C$ consists of one irreducible conic and a line
intersecting tangentially  with the sequence $(2,2)$ or $(2,2,1)$. For $(2,2,1)$, we have two ways to take blow ups as we mentioned right before Lemma~\ref{wt2-1}. For the assertion, we must take blow ups over three points over the line, not the conic.
 \qed
\begin{Lem}\label{wt2-5}
If the divisor $\tilde{D}$ has a $0$-curve with multiplicity $2$ as a
maximal self-intersection number and the $0$-curve intersects only one component of $\tilde{D}$, then  the dual graph is  obtained  by suitable blow ups from the following:
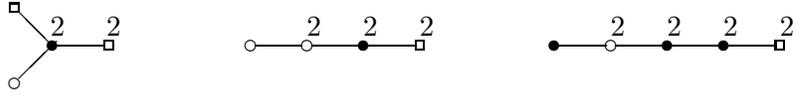
\begin{center}
\begin{minipage}[m]{.2\linewidth}\setlength{\unitlength}{.28mm}\setlength{\unitlength}{.25mm}
\begin{center}

\end{center}
\end{minipage}
\end{center}
\end{Lem}
\pf For the divisor $\tilde{D}$ to have a $0$-curve with
multiplicity $2$, the curve $C$ must consist  of  one double line
and one single line. For such a double line to be a $0$-curve, the sequence must be
$(3,1,1)$, $(3, 1)$ or $(3)$.
 \qed

\begin{Lem}\label{wt2-6}
On a weak del Pezzo surface of degree $4$, there is no effective anticanonical divisor corresponding
to the dual graph
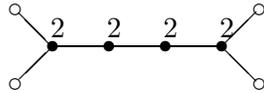
\begin{center}
\begin{minipage}[m]{.28\linewidth}\setlength{\unitlength}{.25mm}
\begin{center}
%
\end{center}

\end{minipage}
\end{center}
\end{Lem}
\pf Contracting five $-1$-curves successively, we get a $0$-curve on
$\mathbb{P}^2$. This is a contradiction.~\qed

\begin{Prop}\label{proposition:wt2}
If the dual graph for an effective anticanonical divisor  on a weak del Pezzo surface  has a vertex of weight $2$ as a maximum,
then it is exactly one of the following:

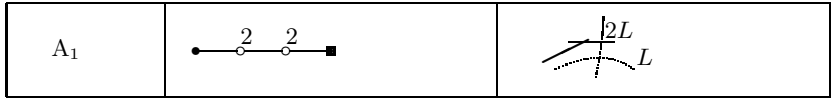
\footnotesize
\begin{center}
\begin{longtable}{|l|c|c|}
\hline \multicolumn{3}{|c|}
{Degree $1$} \\
\hline

\phantom{$\surd$\ $2\textnormal{E}_1+\textnormal{E}_1$}&Configuration& Example \\

\hline  $\surd$\ $\textnormal{D}_8$ &
\begin{minipage}[m]{.25\linewidth}\setlength{\unitlength}{.20mm}
\begin{center}
\begin{picture}(500,60)(90,10)
\put(130,60){\line(1,-1){20}} \put(130,20){\line(1,1){20}}
\put(150,40){\line(1,0){120} }
\put(270,40){\line(1,1){18.3}} \put(270,40){\line(1,-1){18.3}}

\put(290,60){\circle*{5}} \put(290,20){\circle{5}}

\multiput(150,40)(30,0){5}{\circle*{5}} \put(130,60){\circle*{5}}
\multiput(130,20)(100,0){1}{\circle*{5}}


\multiput(149,45)(30,0){5}{$2$}
\end{picture}
\end{center}

\end{minipage}&
\begin{minipage}[m]{.25\linewidth}\setlength{\unitlength}{.20mm}
\begin{center}
\begin{picture}(400,50)(-10,15)

\qbezier[15](124,46)(115,38)(106,30)
\put(128,42){\mbox{\footnotesize $Q$}}
\qbezier[20](18,60)(50,40)(26,30)
\put(12,20){\mbox{\footnotesize $L$}}

\put(118,33){\line(1,0){30}}
\put(10,55){\line(1,0){30}}

\thicklines
\put(0,30){\line(2,1){30}}
\put(20,45){\line(2,-1){30}}
\put(40,30){\line(2,1){30}}
\put(60,45){\line(2,-1){30}}
\put(80,30){\line(2,1){30}} \put(100,45){\line(2,-1){30}}

\end{picture}
\end{center}

\end{minipage}\\
\hline
\end{longtable}
\end{center}

\begin{center}
\begin{longtable}{|l|c|c|}
\hline \multicolumn{3}{|c|}
{Degree $2$} \\
\hline

\phantom{$\surd$\ $2\textnormal{E}_1+\textnormal{E}_1$}&Configuration& Example \\

\hline  $\surd$\ $\textnormal{D}_6+\textnormal{A}_1$ &
\begin{minipage}[m]{.25\linewidth}\setlength{\unitlength}{.20mm}
\begin{center}
\begin{picture}(500,60)(120,10)
\put(130,60){\line(1,-1){20}} \put(130,20){\line(1,1){20}}
\put(150,40){\line(1,0){117.5} }

\multiput(150,40)(30,0){4}{\circle*{5}} \put(130,60){\circle*{5}}
\multiput(130,20)(100,0){1}{\circle*{5}}


\multiput(149,45)(30,0){5}{$2$} \put(270,40){\circle{5}}
\put(272.5,40){\line(1,0){27.5}} \put(300,40){\circle*{5}}

\end{picture}
\end{center}

\end{minipage}
&
\begin{minipage}[m]{.25\linewidth}\setlength{\unitlength}{.20mm}
\begin{center}
\begin{picture}(400,60)(-10,-4)

\qbezier[35](0,40)(55,20)(40,0)
\put(50,0){\mbox{\footnotesize $2L$}}
\qbezier[15](40,45)(55,40)(70,30)
\put(41,46){\mbox{\footnotesize $L$}}

\put(18,12){\line(1,0){30}}
\put(80,33){\line(1,0){30}}

\thicklines
\put(0,30){\line(2,1){30}}
\put(20,45){\line(2,-1){30}}
\put(40,30){\line(2,1){30}}
\put(60,45){\line(2,-1){30}}

\put(0,0){\line(2,1){30}}

\end{picture}

\end{center}
\end{minipage}
\\
\hline
\end{longtable}
\end{center}

\begin{center}
\begin{longtable}{|l|c|c|}
\hline \multicolumn{3}{|c|}
{Degree $3$} \\
\hline

\phantom{$\surd$\ $2\textnormal{E}_1+\textnormal{E}_1$}&Configuration& Example \\

\hline  $\surd$\ $\textnormal{A}_5+\textnormal{A}_1$ &
\begin{minipage}[m]{.25\linewidth}\setlength{\unitlength}{.20mm}
\begin{center}
\begin{picture}(500,60)(115,10)
\put(130,60){\line(1,-1){20}} \put(131.7,21.7){\line(1,1){20}}
\put(150,40){\line(1,0){117.5} }
\put(272.5,40){\line(1,0){30}} \put(300,40){\circle*{5}}
\put(270,40){\circle{5}}

\multiput(150,40)(30,0){4}{\circle*{5}} \put(130,60){\circle*{5}}
\multiput(130,20)(100,0){1}{\circle{5}}

\multiput(149,45)(30,0){5}{$2$}
\end{picture}
\end{center}

\end{minipage}&
\begin{minipage}[m]{.25\linewidth}\setlength{\unitlength}{.20mm}
\begin{center}
\begin{picture}(400,50)(-10,15)
\qbezier[15](45,42)(35,35)(25,25)
\put(41,45){\mbox{\footnotesize $L$}}
\qbezier[15](120,50)(120,40)(120,30)
\put(125,45){\mbox{\footnotesize $Q$}}

\put(98,42){\line(1,0){30}}

\thicklines
\put(0,30){\line(2,1){30}}
\put(20,45){\line(2,-1){30}}
\put(40,30){\line(2,1){30}}
\put(60,45){\line(2,-1){30}}
\put(80,30){\line(2,1){30}}

\end{picture}
\end{center}
\end{minipage}\\
\hline
\end{longtable}
\end{center}

\begin{center}
\begin{longtable}{|l|c|c|}
\hline \multicolumn{3}{|c|}
{Degree $4$} \\
\hline

\phantom{$\surd$\ $2\textnormal{E}_1+\textnormal{E}_1$}&Configuration& Example \\

\hline  $\surd$\ $\textnormal{A}_3+2\textnormal{A}_1$ &
\begin{minipage}[m]{.25\linewidth}\setlength{\unitlength}{.20mm}
\begin{center}
\begin{picture}(500,60)(90,10)

\put(100,40){\line(1,0){27.5} }
\put(132.5,40){\line(1,0){25}} \put(162.5,40){\line(1,0){25}}
\put(280,40){\line(-1,0){27.5}}\put(190,40){\line(1,0){27.5}}\put(220,40){\line(1,0){27.5}}
\multiput(100,40)(180,0){2}{\circle*{5}}
\multiput(160,40)(30,0){3}{\circle*{5}}
\multiput(130,40)(120,0){2}{\circle{5}}
\multiput(130,45)(30,0){5}{$2$}
\end{picture}
\end{center}

\end{minipage}&
\begin{minipage}[m]{.25\linewidth}\setlength{\unitlength}{.20mm}
\begin{center}
\begin{picture}(400,50)(-10,5)

\put(17,20){\line(1,0){30}}
\put(40,32){\line(1,0){30}}

\thicklines
\put(0,30){\line(2,1){30}}
\put(20,45){\line(2,-1){30}}

\put(0,8){\line(2,1){30}}

\thinlines
\qbezier[25](2,39)(40,42)(40,8)
\put(-5,42){\mbox{\footnotesize $2L$}}
\thinlines\qbezier[12](60,20)(60,30)(60,40)
\put(62,35){\mbox{\footnotesize $L$}}

\end{picture}
\end{center}
\end{minipage}\\
\hline
\end{longtable}
\end{center}

\begin{center}
\begin{longtable}{|l|c|c|}
\hline \multicolumn{3}{|c|}
{Degree $5$} \\
\hline

\phantom{$\surd$\ $2\textnormal{E}_1+\textnormal{E}_1$}&Configuration& Example \\

\hline  $\surd$\ $\textnormal{A}_1$ &
\begin{minipage}[m]{.25\linewidth}\setlength{\unitlength}{.20mm}
\begin{center}
\begin{picture}(500,60)(90,10)

\put(100,40){\line(1,0){27.5} }
\put(132.5,40){\line(1,0){25}} \put(162.5,40){\line(1,0){25}}
\multiput(100,40)(90,0){1}{\circle*{5}}
\multiput(130,40)(30,0){2}{\circle{5}}\put(190,40){\circle{5}}
\multiput(130,45)(30,0){2}{$2$}
\end{picture}
\end{center}

\end{minipage}
&
\begin{minipage}[m]{.25\linewidth}\setlength{\unitlength}{.20mm}
\begin{center}
\begin{picture}(400,50)(-20,5)

\put(0,20){\line(1,0){30}}
\put(0,25){\line(1,0){30}}
\put(18,42){\line(1,0){30}}

\thicklines

\put(0,30){\line(2,1){30}}

\thinlines\qbezier[20](40,55)(40,31)(40,15)
\put(43,47){\mbox{\footnotesize $2L$}}
\thinlines\qbezier[28](8,10)(8,40)(50,35)
\put(-5,5){\mbox{\footnotesize $L$}}

\end{picture}
\end{center}
\end{minipage}
\\
\hline
\end{longtable}
\end{center}

\begin{center}
\begin{longtable}{|l|c|c|}
\hline \multicolumn{3}{|c|}
{Degree $6$} \\
\hline

\phantom{$\surd$\ $2\textnormal{E}_1+\textnormal{E}_1$}&Configuration& Example \\
\hline

\phantom{$\surd$}\ $\textnormal{A}_2+\textnormal{A}_1$ &
\begin{minipage}[m]{.25\linewidth}\setlength{\unitlength}{.20mm}
\begin{center}
\begin{picture}(500,60)(120,10)
\put(132,58){\line(1,-1){18}} \put(132,22){\line(1,1){20}}
\put(150,40){\line(1,0){57.5} }
\put(212.5,40){\line(1,0){25}} \put(240,40){\circle*{5}}
\multiput(150,40)(30,0){2}{\circle*{5}} 
\multiput(128,18)(100,0){1}{\framebox(4,4)}\put(210,40){\circle{5}}\multiput(128,58)(100,0){1}{\framebox(4,4)}
\multiput(149,45)(30,0){3}{$2$}
\end{picture}
\end{center}

\end{minipage}&
\begin{minipage}[m]{.25\linewidth}\setlength{\unitlength}{.20mm}
\begin{center}
\begin{picture}(400,50)(-20,10)

\put(38,33){\line(1,0){30}}

\thicklines
\put(0,30){\line(2,1){30}}
\put(20,45){\line(2,-1){30}}

\thinlines\qbezier[15](50,45)(55,32)(60,20)
\put(60,10){\mbox{\footnotesize $L$}}
\thinlines\qbezier[12](10,45)(10,32)(10,20)
\put(10,10){\mbox{\footnotesize $L$}}
\thinlines\qbezier[12](5,45)(5,32)(5,20)
\put(-7,40){\mbox{\footnotesize $L$}}
\end{picture}
\end{center}
\end{minipage}\\
\hline
\end{longtable}
\end{center}

\begin{center}
\begin{longtable}{|l|c|c|}
\hline \multicolumn{3}{|c|}
{Degree $7$} \\
\hline

\phantom{$\surd$\ $2\textnormal{E}_1+\textnormal{E}_1$}&Configuration& Example \\
\hline \phantom{$\surd$}\ $\textnormal{A}_1$ &
\begin{minipage}[m]{.25\linewidth}\setlength{\unitlength}{.20mm}
\begin{center}
\begin{picture}(500,55)(120,12)
 \put(188,38){\framebox(4,4)} \put(130,40){\circle*{5}}
 \put(158,38){\framebox(4,4)}\put(190,40){\circle*{4}}
\put(130,40){\line(1,0){28}} \put(162,40){\line(1,0){25}}
\put(159,45){$2$}
\end{picture}
\end{center}

\end{minipage}
&
\begin{minipage}[m]{.25\linewidth}\setlength{\unitlength}{.20mm}
\begin{center}
\begin{picture}(400,50)(-20,10)

\put(20,43){\line(1,0){30}}

\thicklines
\put(0,30){\line(2,1){30}}

\thinlines\qbezier[25](0,20)(30,40)(35,15)
\put(37,20){\mbox{\footnotesize $L$}}
\thinlines\qbezier[15](8,15)(8,27)(8,55)
\put(-12,42){\mbox{\footnotesize $2L$}}

\end{picture}
\end{center}
\end{minipage}
\\
\hline
\end{longtable}
\end{center}
\normalsize
\end{Prop}
\pf
Starting from the dual graphs for the degree $1$ case,  we apply
$\textbf{P-1}$ and $\textbf{P-2}$
successively to get the possible dual graphs. Then from the obtained possible dual graphs, we exclude the dual graphs that violate the properties in Lemmas~\ref{wt2-1},~\ref{wt2-2},~\ref{wt2-3},~\ref{wt2-4},~\ref{wt2-5} and~\ref{wt2-6}.
\qed

\begin{Prop}
If the dual graph for an effective anticanonical divisor on a weak del Pezzo surface  has a vertex of weight $1$ as a maximum, the dual graph is circular except the cases where the effective anticanonical divisor consists of either three curves intersecting transversally at a single point or  two curves intersecting tangentially with intersection number $2$ at a single point.
\end{Prop}
\pf The assertion holds for a weak del Pezzo surface of degree $1$ (see \cite{P01}). Notice that we can apply only \textbf{P-1}. Then the assertion is clear.
\qed

\begin{Prop}\label{rational} Let $S$ be a del Pezzo surface of degree $d\geq 2$.
\begin{enumerate}
\item If the surface $S$ has only one singular point that is of type $\textnormal{A}_1$,
then
 there is an
    effective anticanonical divisor on its minimal resolution
    consisting of one $(-3+d)$-curve, one $-1$-curve, and one
    $-2$-curve intersecting transversally at a single
    point.

\item If the surface $S$  has only one singular point that is
    of type $\textnormal{A}_2$, then there is an
    effective anticanonical divisor on its minimal resolution
    consisting of one $(-2+d)$-curve and two $-2$-curves
    intersecting transversally at a single point.
\item If the surface $S$ is of degree $2$ and it has only two
    singular points that are of type $\textnormal{A}_1$ or
    $\textnormal{A}_2$, then there is at least one of the
    following effective anticanonical divisors:
\begin{itemize}
\item One $-1$-curve, one $-1$-curve, and one $-2$-curve
    intersecting transversally at a single point;
\item Two $-2$-curves and one $0$-curve intersecting
    transversally at a single point.
\end{itemize}
\end{enumerate}
\end{Prop}
\pf
Considering  successive suitable blow ups of a cuspidal cubic, three
lines intersecting at a single point, and a conic and a line intersecting
tangentially on $\mathbb{P}^2$, we can easily obtain the assertions.
\qed

\section{Log canonical threshold.}

In this section, we prove Theorem~\ref{theorem:main}. For given singularity types, we  consider all the possible effective anticanonical divisors. However, we do not have to consider all effective anticanonical divisors. It turns out that we have only to consider those which appear in the tables of Propositions~\ref{proposition:wt6}, \ref{proposition:wt4}, \ref{proposition:wt3}, \ref{proposition:wt2} with the mark $\surd$ and those described in Proposition~\ref{rational}. In what follows we explain the reason.

\begin{Prop}\label{change}
Let $D$ be
 an effective
 anticanonical divisor  that contains a $-2$-curve
 on a weak del Pezzo surface $\tilde{S}$.
Write $D=D_1+D_2$ where $D_1$ and $D_2$ are effective divisors
with $D_2^2\geq 0$ and $D_1\cdot D_2=2$. Suppose that there is a $-1$-curve $L$
 on $\tilde{S}$ such that $D^2_2\geq D_2\cdot L \geq 0$. Then there is an effective divisor $D_3$
  such that $D_2$ is linearly equivalent to $L+D_3$. In
  particular,~$D_3^2\leq D_2^2-1$.
  In addition, if $D^2_2=0$, then the divisor $D_3$ consists only of negative curves.
\end{Prop}
\pf By Riemann-Roch,
$$h^0(\tilde{S},\mathcal{O}_{\tilde{S}}(D_2))\geq \frac{D_1\cdot
D_2}{2}+D_2^2+1=2+D_2^2,$$ since
$h^0(\tilde{S},\mathcal{O}(K_{\tilde{S}}-D_2))=0$. Consider the
exact sequence $$0 \rightarrow
\mathcal{O}_{\tilde{S}}(D_2-L)\rightarrow\mathcal{O}_{\tilde{S}}(D_2)\rightarrow
\mathcal{O}_L(D_2|_L)\rightarrow 0.$$ We then see that
$D_2-L$ is linearly equivalent to an effective divisor $D_3$ since
$$h^0(\tilde{S},\mathcal{O}_{\tilde{S}}(D_2-L))\geq h^0(\tilde{S},\mathcal{O}_{\tilde{S}}(D_2))-h^0(L,\mathcal{O}_L(D_2|_L))\geq 1+D_2^2-D_2\cdot L\geq 1.$$
Furthermore, $D_3^2=(D_2-L)^2=D^2_2+L^2-2D_2\cdot L\leq D^2_2-1$.

If $D^2_2=0$, then $-K_{\tilde{S}}\cdot D_3=(D_1+D_2)\cdot (D_2-L)=1$.
Since $D_3$ is a part of an effective anticanonical divisor with a $-2$-curve, it consists of rational curves. Therefore, we can conclude that $D_3$  consists of one $-1$-curve and $-2$-curves.
 \qed

\begin{Lem}\label{line} Let $\tilde{S}$ be a weak del Pezzo surface of degree $\leq 7$. Then every non-negative nonsingular rational curve $C$ that appears in an effective anticanonical divisor on $\tilde{S}$ has a $-1$-curve $L$ with $L\cdot C=0$.
\end{Lem}
\pf We have a sequence of blow ups $\pi:\tilde{S}\to\mathbb{P}^2$.  Then the curve $C$ is the strict transform of
a line,  a conic or a singular cubic on $\mathbb{P}^2$. If the curve $C$ comes from a singular cubic, then the singular point must be the center of an exceptional divisor of $\pi$. Suppose
that there is  at least two $-1$-curves on $\tilde{S}$ that are exceptional curves of $\pi$ and $C$ intersects with all the $-1$-curves. Choose two points among the centers of all exceptional curves of $\pi$ in such a way that if the curve $C$ comes from a singular cubic, the singular point is one of the two points.   The strict
transform of the line passing through the chosen two points does not meet $C$ and it is a
$-1$-curve; otherwise the curve $C$ would not be irreducible. On the other hand, if there is  exactly one  exceptional $-1$-curve of $\pi$ such that
$C$ intersects with this line, then the center of the exceptional
curves on $\mathbb{P}^2$ is one point  and the curve $C$
must be the strict transform of either a conic or a singular cubic. If $C$ is the strict transform of a conic, consider the line on
$\mathbb{P}^2$  that is tangent to a conic at the point of the center. If $C$ is the strict transform of a singular cubic, consider the line on
$\mathbb{P}^2$ that is one of components of the tangent cone of the singular point of the cubic. Then the strict
transform of the tangent line on $\tilde{S}$  does not intersect  $C$
and it is  a $-1$-curve on $\tilde{S}$.\qed

From now on, for a divisor $D$ on $\tilde{S}$  we define
$$\mult(D):=\max\{\mult_C(D) : C \textnormal{ is an irreducible curve on } \tilde{S}\}.$$

\begin{Prop}\label{negative-curves-only}

An effective anticanonical divisor $D$ containing a $-2$-curve on $\tilde{S}$ is linearly equivalent to an effective anticanonical divisor $D'$ containing a $-2$-curve such that
\begin{itemize}
\item it consists of only $-1$-curves and $-2$-curves;
\item $\mult(D')\geq \mult(D)$.
\end{itemize}
\end{Prop}
\pf Suppose that the divisor $D$ contains a non-negative curve $C$. Let $m=\mult_C(D)$. By Lemma~\ref{line}, there is a $-1$-curve $L$ with $C\cdot L=0$. Then we obtain an effective divisor $F$ with~$F^2\leq C^2-1$ such that $C$ is linearly equivalent to $L+F$
by Lemma~\ref{change}. We replace $D$ by the effective divisor $D-mC+m(L+F)$. Note that $\mult(D)\leq \mult(D-mC+m(L+F))$. Repeating this procedure finitely many times we get the required effective anticanonical divisor.~\qed

\begin{Lem}\label{line2} Let $D$ be an effective anticanonical divisor on a weak del Pezzo surface $\tilde{S}$. Write~$D=D_1+D_2$, where $D_1$ and $D_2$ are effective divisors.
Suppose that the divisor $D_2$ has the dual graph as follows:\\
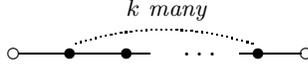
\begin{center}
\hspace{50mm}
\begin{minipage}[m]{.2\linewidth}\setlength{\unitlength}{.28mm}\setlength{\unitlength}{.25mm}
\begin{center}
\begin{picture}(100,0)(260,40)
\put(102.5,40){\line(1,0){25} }
\put(132.5,40){\line(1,0){40}} \put(190,39){$\ldots$}
\multiput(130,40)(30,0){2}{\circle*{5}}\put(220,40){\line(1,0){33}}
\multiput(230,40)(30,0){1}{\circle*{5}}
\multiput(100,40)(60,0){1}{\circle{5}}\put(255.5,40){\circle{5}}
\thinlines\qbezier[35](133,45)(180,60)(227,45)
\put(160,60){\mbox{\footnotesize $k$ many}}
\end{picture}
\end{center}
\end{minipage}
\end{center}
where $k\geq 0$.
  In addition, we suppose that the divisor $D_1$ has a $-2$-curve that is not a component of $D_2$ and  does not intersect at least one $-1$-curve in $D_2$. Then there exists a $-1$-curve
  that does not intersect the divisor $D_2$.
\end{Lem}
\pf
Let $E$ be a $-2$-curve in $D_1$ that is not a component of $D_2$ and  does not intersect at least one $-1$-curve in $D_2$.
By blowing down $-1$-curves $k+1$ times from the $-1$-curve that does not intersect $E$,
we get a $0$-curve that is contained in an effective anticanonical divisor containing $-2$-curve on a new weak del Pezzo surface.

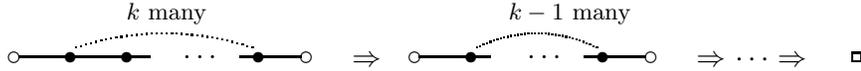
\begin{center}
\hspace{50mm}
\begin{minipage}[m]{.2\linewidth}\setlength{\unitlength}{.28mm}\setlength{\unitlength}{.25mm}
\begin{center}
\begin{picture}(100,0)(260,55)
\put(102.5,40){\line(1,0){25} }
\put(132.5,40){\line(1,0){40}} \put(190,39){$\ldots$}
\multiput(130,40)(30,0){2}{\circle*{5}}\put(220,40){\line(1,0){33}}
\multiput(230,40)(30,0){1}{\circle*{5}}
\multiput(100,40)(60,0){1}{\circle{5}}\put(255.5,40){\circle{5}}
\put(280,35){$\Rightarrow$}
\qbezier[35](133,45)(180,60)(227,45)
\put(160,60){\mbox{\footnotesize $k$ many}}
\end{picture}
\end{center}
\end{minipage}
\begin{minipage}[m]{.2\linewidth}\setlength{\unitlength}{.28mm}\setlength{\unitlength}{.25mm}
\begin{center}
\begin{picture}(0,0)(260,55)

\put(132.5,40){\line(1,0){38}} \put(190,39){$\ldots$}
\multiput(160,40)(30,0){1}{\circle*{5}}\put(220,40){\line(1,0){33}}
\multiput(230,40)(30,0){1}{\circle*{5}}
\multiput(130,40)(60,0){1}{\circle{5}}
\put(255.5,40){\circle{5}}
\put(280,35){$\Rightarrow$ $\cdots$}
\thinlines\qbezier[35](163,45)(195,60)(227,45)
\put(180,60){\mbox{\footnotesize $k-1$ many}}
\end{picture}
\end{center}
\end{minipage}
\begin{minipage}[m]{.2\linewidth}\setlength{\unitlength}{.28mm}\setlength{\unitlength}{.25mm}
\begin{center}
\begin{picture}(0,0)(80,55)
\put(10,35){$\Rightarrow$ }
\put(50,38){\framebox(4,4)}
\end{picture}
\end{center}
\end{minipage}
\end{center}\medskip
We then apply Lemma~\ref{line} to obtain the required $-1$-curve.
\qed

\begin{Lem}\label{line3} Let $D$ be an effective anticanonical divisor
 on a weak del Pezzo surface $\tilde{S}$
of degree~$2$. Write $D=D_1+D_2$ where $D_1$ and $D_2$ are
effective divisors.
Suppose that the divisor $D_2$ has the dual graph as follows:\\
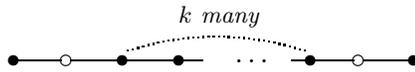
\begin{center}
\hspace{50mm}
\begin{minipage}[m]{.2\linewidth}\setlength{\unitlength}{.28mm}\setlength{\unitlength}{.25mm}
\begin{center}
\begin{picture}(72,0)(260,40)
\put(72.5,40){\line(1,0){25} }
\multiput(72,40)(30,0){1}{\circle*{5}}
\multiput(285,40)(30,0){1}{\circle*{5}}
\put(102.5,40){\line(1,0){25} }
\put(132.5,40){\line(1,0){40}} \put(190,39){$\ldots$}
\multiput(130,40)(30,0){2}{\circle*{5}}\put(220,40){\line(1,0){33}}
\put(258,40){\line(1,0){25}}
\multiput(230,40)(30,0){1}{\circle*{5}}
\multiput(100,40)(60,0){1}{\circle{5}}\put(255.5,40){\circle{5}}
\thinlines\qbezier[35](133,45)(180,60)(227,45)
\put(160,60){\mbox{\footnotesize $k$ many}}
\end{picture}
\end{center}
\end{minipage}
\end{center}
where $k\geq 0$.
  In addition, we suppose that the divisor $D_1$ has a $-2$-curve that is not a component of $D_2$ and  not connected to any $-2$-curve between two $-1$-curves in $D_2$. Then there exist a $-1$-curve
  that does not intersect the divisor $D_2$.
\end{Lem}
\pf
Let $E$ be a $-2$-curve in $D_1$ that is not a component of $D_2$ and  not connected to any $-2$-curve between two $-1$-curves in $D_2$. Applying Lemma~\ref{line2} to the divisor obtained by subtracting two $-2$-curves at the ends of $D_2$ from $D_2$, we see that there is a $-1$-curve $L$  that intersects either $E$ or
one of two $-2$-curves at the ends of $D_2$. If the $-1$-curve $L$ intersects $E$,
then we are done. Therefore, we suppose that the curve $L$ intersects
one of two $-2$-curves at the ends of $D_2$.

If $k=0$, then we contract $L$ and then blow down $-1$-curves as
follows:

\begin{center}
\hspace{50mm}
\begin{minipage}[m]{.2\linewidth}\setlength{\unitlength}{.28mm}\setlength{\unitlength}{.25mm}
\begin{center}
\begin{picture}(72,0)(180,40)
\put(72.5,40){\line(1,0){25} }
\multiput(72,40)(30,0){1}{\circle*{5}}
\put(102.5,40){\line(1,0){25} }
\put(132.5,40){\line(1,0){25}}
\multiput(130,40)(30,0){2}{\circle{5}}
\multiput(100,40)(60,0){1}{\circle{5}} \put(180,36){$\Rightarrow$}
\end{picture}
\end{center}
\end{minipage}
\begin{minipage}[m]{.2\linewidth}\setlength{\unitlength}{.28mm}\setlength{\unitlength}{.25mm}
\begin{center}
\begin{picture}(100,0)(180,40)
\put(102.5,40){\line(1,0){25} }
\put(132.5,40){\line(1,0){25}}
\multiput(160,40)(30,0){1}{\circle{5}}
\put(128,38){\framebox(4,4)}
\multiput(100,40)(60,0){1}{\circle{5}} \put(180,36){$\Rightarrow$}
\end{picture}
\end{center}
\end{minipage}
\begin{minipage}[m]{.1\linewidth}\setlength{\unitlength}{.28mm}\setlength{\unitlength}{.25mm}
\begin{center}
\begin{picture}(0,0)(120,40)
\put(4.5,40){\line(1,0){25}} \multiput(32,40)(30,0){1}{\circle{5}}
\put(0,38){\framebox(4,4)} \multiput(2,40)(60,0){1}{\circle*{4}}
\put(52,36){$\Rightarrow$}
\end{picture}
\end{center}
\end{minipage}
\begin{minipage}[m]{.1\linewidth}\setlength{\unitlength}{.28mm}\setlength{\unitlength}{.25mm}
\begin{center}
\begin{picture}(0,0)(100,40)
\put(0,38){\framebox(4,4)} \put(2,40){\circle{4}}
\end{picture}
\end{center}
\end{minipage}
\end{center}
If $k=1$, then we blow down $-1$-curves as follows:
\begin{center}
\hspace{30mm}
\begin{minipage}[m]{.2\linewidth}\setlength{\unitlength}{.28mm}\setlength{\unitlength}{.25mm}
\begin{center}
\begin{picture}(72,0)(200,40)
\put(72.5,40){\line(1,0){25} }
\multiput(72,40)(30,0){1}{\circle*{5}}
\put(102.5,40){\line(1,0){25} }
\put(132.5,40){\line(1,0){25}} \put(162.5,40){\line(1,0){25}}
\multiput(130,40)(30,0){1}{\circle*{5}}
\multiput(100,40)(60,0){2}{\circle{5}} \put(190.5,40){\circle*{5}}
\put(210,36){$\Rightarrow$}
\end{picture}
\end{center}
\end{minipage}
\begin{minipage}[m]{.2\linewidth}\setlength{\unitlength}{.28mm}\setlength{\unitlength}{.25mm}
\begin{center}
\begin{picture}(72,0)(160,40)
\put(72.5,40){\line(1,0){25} }
\multiput(72,40)(30,0){1}{\circle*{5}}
\put(102.5,40){\line(1,0){25} }
\put(132.5,40){\line(1,0){25}}
\multiput(130,40)(30,0){2}{\circle{5}}
\multiput(100,40)(60,0){1}{\circle{5}} \put(180,36){$\Rightarrow$}
\end{picture}
\end{center}
\end{minipage}
\begin{minipage}[m]{.15\linewidth}
\setlength{\unitlength}{.28mm}\setlength{\unitlength}{.25mm}
\begin{center}
\begin{picture}(0,0)(100,40)
\put(2.5,40){\line(1,0){25} }
\put(32.5,40){\line(1,0){25}}
\multiput(60,40)(30,0){1}{\circle{5}} \put(28,38){\framebox(4,4)}
\multiput(0,40)(60,0){1}{\circle{5}} \put(80,36){$\Rightarrow$}
\end{picture}
\end{center}
\end{minipage}
\begin{minipage}[m]{.1\linewidth}\setlength{\unitlength}{.28mm}\setlength{\unitlength}{.25mm}
\begin{center}
\begin{picture}(0,0)(70,40)
\put(4.5,40){\line(1,0){25}} \multiput(32,40)(30,0){1}{\circle{5}}
\put(0,38){\framebox(4,4)} \multiput(2,40)(60,0){1}{\circle*{4}}
\put(52,36){$\Rightarrow$}
\end{picture}
\end{center}
\end{minipage}
\begin{minipage}[m]{.1\linewidth}\setlength{\unitlength}{.28mm}\setlength{\unitlength}{.25mm}
\begin{center}
\begin{picture}(0,0)(50,40)
\put(0,38){\framebox(4,4)} \put(2,40){\circle{4}}
\end{picture}
\end{center}
\end{minipage}
\end{center}
In both cases, we get a $2$-curve that is contained in an
effective anticanonical divisor on a new weak del Pezzo surface.
The weak del-Pezzo surface is of degree $\leq 7$, Lemma~\ref{line}
implies the assertion.

Now we suppose $k\geq 2$. By contracting  $L$ we get a divisor on
a new weak del Pezzo surface whose dual graph is  as follows:
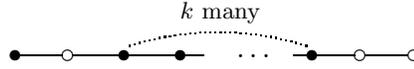
\begin{center}
\hspace{50mm}
\begin{minipage}[m]{.2\linewidth}\setlength{\unitlength}{.28mm}\setlength{\unitlength}{.25mm}
\begin{center}
\begin{picture}(72,0)(260,47)
\put(72.5,40){\line(1,0){25} }
\multiput(72,40)(30,0){1}{\circle*{5}}
\multiput(284,40)(30,0){1}{\circle{5}}
\put(102.5,40){\line(1,0){25} }
\put(132.5,40){\line(1,0){40}} \put(190,39){$\ldots$}
\multiput(130,40)(30,0){2}{\circle*{5}}\put(220,40){\line(1,0){33}}
\put(258,40){\line(1,0){23}}
\multiput(230,40)(30,0){1}{\circle*{5}}
\multiput(100,40)(60,0){1}{\circle{5}}\put(255.5,40){\circle{5}}
\thinlines\qbezier[35](133,45)(180,60)(227,45)
\put(160,60){\mbox{\footnotesize $k$ many}}
\end{picture}
\end{center}
\end{minipage}
\end{center}
Then we blow down $-1$-curves $k$ times from the $-1$-curve on the left to the right. Then we obtain a divisor on a new weak del Pezzo surface whose dual graph is  as follows:
\begin{center}
\hspace{50mm}
\begin{minipage}[m]{.2\linewidth}\setlength{\unitlength}{.28mm}\setlength{\unitlength}{.25mm}
\begin{center}
\begin{picture}(121,0)(260,37)
\
\put(132.5,40){\line(1,0){25}}\put(162.5,40){\line(1,0){25}}\put(192.5,40){\line(1,0){25}}
\multiput(160,40)(30,0){3}{\circle{5}}\multiput(121,36)(30,0){1}{\mbox{
$\diamond$}}

\end{picture}
\end{center}
\end{minipage}
\end{center}
where $\diamond$ is the curve from the $-2$-curve at the left end of $D_2$. Again we contract the $-1$-curve in the middle. Then we get  a divisor on a new weak del Pezzo surface whose dual graph is  as follows:
\begin{center}
\begin{minipage}[m]{.2\linewidth}\setlength{\unitlength}{.28mm}\setlength{\unitlength}{.25mm}
\begin{center}
\begin{picture}(-100,0)(260,37)
\ \put(132.5,40){\line(1,0){25}}\put(162.5,40){\line(1,0){25}}
\multiput(158,38)(30,0){2}{\framebox(4,4)}\multiput(121,36)(30,0){1}{\mbox{
$\diamond$}}
\end{picture}
\end{center}
\end{minipage}

\end{center}

Suppose that all $-1$-curve on this surface are connected to this divisor. Then we contract all the $-1$-curves so that we obtain either $\mathbb{P}^2$ or $\mathbb{P}^1\times\mathbb{P}^1$. However, since the $-2$-curves $E$ has never been touched by any $-1$-curve, this is a contradiction. Therefore, there is a $-1$-curve that does not intersect the divisor $D_2$.
 \qed

\begin{Prop}\label{no-1} Let $\tilde{S}$ be a weak del Pezzo surface of degree $d$.
Suppose that there is an effective anticanonical divisor $D$ with $\mult(D)=1$ such that it contains at least one $-2$-curve, consists of only $-1$-curves and $-2$-curves and  satisfies the following conditions:
\begin{itemize}
\item if $d\geq 3$, then it contains at least five curves;
\item if $d=2$, then it contains a chain of at least three $-2$-curves.
\end{itemize}
Then there is an effective anticanonical divisor $D'$ on $\tilde{S}$ with $\mult(D')\geq 2$.
\end{Prop}
\pf
The dual graph of the divisor $D$ must be circular.

First, we suppose that the divisor $D$ contains at least three $-1$-curves.

If either its dual graph contains at least two chains of $-1$-curves or  it contains at least four $-1$-curves and only one chain of $-1$-curves, then we can obtain a divisor $D_2$ from $D$ that satisfies the conditions
of Lemma~\ref{line2}. Furthermore, we can pick a $-1$-curve $L$ from the divisor $D$ not intersecting $D_2$.
Then Proposition~\ref{change} implies the assertion.

If the dual graph contains exactly three $-1$-curves and only one chain of $-1$-curves, we let $L_1$, $L_2$, $L_3$ be the three $-1$-curves with $L_1\cdot L_3=0$. Then the divisor $D_2=L_1+L_2$ must have a $-1$-curve $L$ that dose not intersect $D_2$ by Lemma~\ref{line2}. If the $-1$-curve $L$ intersects the $-2$-curve that intersects $L_3$, then the divisor $L_1+L_2$ must be linearly equivalent to an divisor~$L+R$, where~$R$ is an effective divisor,  by  Proposition~\ref{change}. Then  $D-L_1-L_2+L+R$ is an effective anticanonical divisor whose dual graph has a fork. Therefore, it must have a multiple component. If the curve $L$ intersects the $-2$-curve that intersects $L_1$, then the divisor~$L_2+L_3$ and the $-1$-curve $L$ work for the assertion in the same manner.

Secondly, we suppose that the divisor $D$ contains only two
$-1$-curves. Consider the divisor~$D_1<D$ that consists of the two
$-1$-curves and all the $-2$-curves contained in the chain from
one $-1$-curve to the other not longer than the other side. Note
that $D_1$ contains no $-2$-curve if two $-1$-curves are
connected. Lemma~\ref{line3} shows that there is a $-1$-curve $L$
that intersects a $-2$-curve not in $D_1$ and not connected to the
$-1$-curves in $D_1$. Proposition~\ref{change} then implies that
the divisor $D$ is linearly equivalent to an effective
anticanonical divisor whose dual graph has a fork. This completes
the proof. \qed

Therefore, Propositions~\ref{negative-curves-only} and~\ref{no-1} show that we do not have to consider effective anticanonical divisors without a multiple curve except those in Proposition~\ref{rational}. Furthermore, applying Proposition~\ref{change} with Lemmas~\ref{line} and~\ref{line2} to effective anticanonical divisors with a multiple curve
we are able to obtain a short list of effective anticanonical divisors with a multiple curve to be considered for the first log canonical thresholds. Such divisors are marked by $\surd$ in the tables. This short list gives a proof of Theorem~\ref{theorem:main}.

 \footnotesize
\bibliographystyle{amsplain}

\begin{thebibliography}{10}


\bibitem{Ch07a}
I.\,Cheltsov, \emph{Log canonical thresholds of del Pezzo surfaces},
Geom. Funct. Anal. \textbf{11} (2008) 1118--1144.
\bibitem{Ch07b}
I.\,Cheltsov, \emph{On singular cubic surfaces}, Asian J. of Math. \textbf{13} (2009) 191--214.



\bibitem{CSh08}
I.\,Cheltsov, K.\,Shramov, \emph{Log canonical thresholds of
smooth Fano threefolds},
Russian Math. Surveys \textbf{63} (2008), 859--958. %





\bibitem{DeKo01}
J.-P.\,Demailly, J.\,Koll\'ar, \emph{Semi-continuity of complex singularity exponents and K\"ahler-Einstein metrics on Fano orbifolds},
 Ann. Sci. \'Ecole Norm. Sup. \textbf{34} (2001), 525--556.%



\bibitem{Dema80}M.~Demazure,
\emph{Surfaces de del {P}ezzo}, S\'eminaire sur les singularit\'es
des surfaces, Lect. Notes in Math., vol. \textbf{777}, Springer-Verlag,
1980,  pp.~23--69.



\bibitem{Furu86}M.~Furushima,
\emph{Singular del {P}ezzo surfaces and analytic
compactifications of $3$-dimensional complex affine space
$\mathbb{C}^3$}, Nagoya Math. J. \textbf{104} (1986), 1--28.

\bibitem{HiWa81}F.~Hidaka, K.~Watanabe,
\emph{Normal {G}orenstein surfaces with ample  anti-canonical
divisor}, Tokyo J. Math. \textbf{4} (1981), no.~2, 319--330.




\bibitem{MiZha88}M.~Miyanishi, D.-Q. Zhang,
\emph{Gorenstein log del {P}ezzo surfaces of rank  one}, J.
Algebra \textbf{118} (1988), no.~1, 63--84.

\bibitem{MiZha93}\bysame,
\emph{Gorenstein log del {P}ezzo surfaces. {I}{I}}, J. Algebra
\textbf{156} (1993), no.~1, 183--193.



\bibitem{P99}J.~Park,
\emph{Birational maps of del {Pezzo} fibrations}, J. Reine
Angew.Math. \textbf{538} (2001), 213--221.

\bibitem{P01}\bysame,
\emph{A note on del {Pezzo} fibrations of degree $1$},
 Comm. Algebra \textbf{31} (2003), no.~12
5755--5768.

\bibitem{Pi77}H.~C. Pinkham, \emph{Simple elliptic singularities},
Del Pezzo Surfaces and Cremona Transformations, Proc Symp. in Pure
Math., \textbf{30} (Several Complex Variables), 1977, pp.69--70.


\bibitem{Pu05}
A.~V. Pukhlikov, \emph{Birational geometry of Fano direct
products}, Izv. Math. \textbf{69} (2005), 1225--1255.




\bibitem{U83}T.~Urabe,
\emph{On singularities on degenerate del {P}ezzo surfaces of
degree  $1,$ $2$}, Singularities, Part 2 (Arcata, Calif., 1981),
Amer. Math. Soc.,  Providence, R.I., 1983, pp.~587--591.



\bibitem{Zha88}D.-Q. Zhang,
\emph{Logarithmic del {P}ezzo surfaces of rank one with
contractible  boundaries}, Osaka J. Math. \textbf{25} (1988),
no.~2, 461--497.

\end{thebibliography}
\providecommand{\bysame}{\leavevmode\hbox to3em{\hrulefill}
\thinspace}

\end{document}